\documentclass[11pt,english]{article}
\usepackage[T1]{fontenc}

\usepackage{graphicx, subfigure}
\expandafter\def\csname ver@subfig.sty\endcsname{}
\usepackage[latin9]{inputenc}
\usepackage{geometry}
\usepackage{amsmath}
\usepackage{appendix}
\usepackage{amssymb}
\usepackage{pifont}
\usepackage{esint}
\usepackage{makecell,multirow,diagbox}
\usepackage{mathtools}
\usepackage{epstopdf}
\usepackage[latin9]{inputenc}
\usepackage{pict2e}
\usepackage{xcolor}%定义了一些颜色
\usepackage{colortbl,booktabs}%第二个包定义了几个*rule
\usepackage{stmaryrd}
\usepackage{esint}
\usepackage[all,cmtip]{xy}
\usepackage{mathtools}
\usepackage{float}
\usepackage{setspace}
\usepackage{authblk}
\usepackage{amsthm}
\usepackage{mathrsfs}
\usepackage{stmaryrd}
\usepackage{bbold}
\usepackage{color}
\usepackage{tikz}
\usepackage[all,cmtip]{xy}
\usepackage{algorithm}  
\usepackage{algorithmic}

\theoremstyle{plain}
\theoremstyle{plain}
\newtheorem{remark}{Remark}
\usepackage{color}
\definecolor{marin}{rgb} {0., 0.3, 0.7}
\definecolor{rouge}{rgb} {0.8, 0., 0.}
\definecolor{sepia}{rgb} {0.8, 0.5, 0.}
\usepackage[colorlinks,citecolor=sepia,linkcolor=marin,bookmarksopen,bookmarksnumbered]{hyperref}

\theoremstyle{definition}

\DeclareSymbolFont{largesymbol}{OMX}{yhex}{m}{n}
\DeclareMathAccent{\Widehat}{\mathord}{largesymbol}{"62}

\makeatletter

\usepackage{babel}
\usepackage{subfig}
\graphicspath{{Dirac/}}
           
\def\Eb{\mathbf{E}}
\def\Fb{\mathbf{F}}
\def\Bb{\mathbf{B}}

\def\Jb{\mathbf{J}}

\def\xb{\mathbf{x}}
\def\vb{\mathbf{v}}
\def\ab{\mathbf{a}}
\def\bb{\mathbf{b}}

\def\etab{\boldsymbol \eta}

\def\te{\text{e}}
\def\tA{\text{A}}
\def\tB{\text{B}}

\def\cH{\mathcal{H}}

\def\be{\begin{equation}}
\def\ee{\end{equation}}
\def\tn{\textnormal}
           
\newcommand{\del}[2]{\frac{ \delta #1}{ \delta #2}}
\newcommand{\pder}[2]{\frac{ \partial #1}{ \partial #2}}

\makeatother
\DeclareMathSizes{14}{14}{9.8}{7}
\usepackage{lipsum}

\begin{document}

\title{ Geometric Particle-In-Cell discretizations of a plasma hybrid model with kinetic ions and mass-less fluid electrons}
\date{}
\author[1]{Yingzhe Li \thanks{Corresponding author: yingzhe.li@ipp.mpg.de}}
\author[1]{Martin Campos Pinto}
\author[1,2]{Florian Holderied}
\author[1]{Stefan Possanner}
\author[1,3]{Eric Sonnendr\"ucker}
\affil[1]{Max Planck Institute for Plasma Physics, Boltzmannstrasse 2, 85748 Garching, Germany}
\affil[2]{Technical University of Munich, Department of Physics, James-Franck-Strasse 1, 85748 Garching, Germany}
\affil[3]{Technical University of Munich, Department of Mathematics, Boltzmannstrasse 3, 85748 Garching, Germany}
\maketitle
%##########################
\begin{abstract}
We explore the possibilities of applying structure-preserving numerical methods to a plasma hybrid model with kinetic ions and mass-less fluid electrons satisfying the quasi-neutrality relation. The numerical schemes are derived by finite element methods in the framework of  finite element exterior calculus
(FEEC) for field variables, particle-in-cell (PIC) methods for the Vlasov equation, and splitting methods in time based on an anti-symmetric bracket proposed. Conservation properties of energy, quasi-neutrality relation, positivity of density, and divergence-free property of the magnetic field are given irrespective of the used resolution and metric. Local quasi-interpolation is used for dealing with the current terms in order to make the proposed methods more efficient. The implementation has been done in the framework of the Python package
STRUPHY~\cite{3}, and has been verified by extensive numerical experiments.

\end{abstract}
\setcounter{tocdepth}{1} %to only see sections

%\tableofcontents

\section{Introduction}

Kinetic models and fluid models are used extensively in computational plasma physics. For kinetic models, classical methods include Particle-In-Cell methods~\cite{pic1, pic2}, Eulerian methods~\cite{Califano}, and semi-Lagrangian methods~\cite{semi1}. Kinetic models with high dimensional phase-spaces are able to include physical processes far from thermal equilibrium, but are computationally expensive. Fluid models such as Magnetohydrodynamics (MHD), on the other hand, have lower dimensions and are more efficient for plasmas with nearly Gaussian velocity distributions.
Hybrid models combine the advantages of kinetic and fluid approaches in situations where parts of the plasma are thermalized (fluid, many collisions) and other parts remain far from thermal equilibrium (kinetic species, few collisions). In this paper, we consider a plasma hybrid model in which ions are treated kinetically and mass-less electrons are described using fluid equations with an extended Ohm's law~\cite{current2D} (including electron pressure gradient and Hall term): 
\begin{equation} \label{model}
\begin{aligned}
\text{kinetic ions:}\qquad &\frac{\partial f}{\partial t} + {\mathbf v} \cdot \frac{\partial f}{\partial \xb} + \frac{q}{m} ({\mathbf E} + {\mathbf v} \times {\Bb}) \cdot \frac{\partial f}{\partial {\mathbf v}} = 0\,,
\\[2mm]
\text{Faraday's law:}\qquad &\frac{\partial \Bb}{\partial t} = - \nabla \times {\mathbf E}\,,\qquad \nabla \cdot \Bb = 0\,,
\\[1mm]
\text{mass-less fluid  electrons:} \qquad & \frac{\partial n_\te}{\partial t} +  \nabla \cdot \left(n_\te \mathbf u  \right) = 0\,,
\\[1mm]
\text{Ohm's law:}  \qquad & {\mathbf E} = - \frac{k_\tB T_\te}{e }\frac{\nabla n_\te^{\gamma}}{n_\te} - \left(\mathbf u - \frac{\mathbf J}{n_\te e}  \right) \times {\Bb} \,.
\end{aligned}
\end{equation}
Here,  
 $\Jb = \frac{1}{\mu_0} \nabla \times {\Bb}$ denotes the plasma current,  $f$ denotes the ion distribution function, quasi-neutrality number density is $n_\te = n$, where $n = \frac{q}{e} \int f\,\tn{d}{\mathbf v} $, ${\mathbf u} = q \int {\mathbf v} f \mathrm{d}{\mathbf v}/en$ is the ion current carrying drift velocity, $-e$ is the electron charge, $q/m$ is the ion charge-to-mass ratio,
$\Eb$ and $\Bb$ stand for the electromagnetic fields, $T_\te = const$ is the fixed electron temperature and $k_\tB$ stands for the Boltzmann constant. The choice $\gamma=1 \ \text{or} \ \gamma = 5/3$ results in an isothermal or adiabatic electron model. In this work, we focus on developing numerical methods for the case  $\gamma=1$, which could apply to the case $\gamma=5/3$ directly. 

By straightforward calculations it is verified that system \eqref{model} conserves the energy when $\gamma = 1$,
\be \label{H}
\mathcal{H} = \underbrace{\frac{m}{2} \int |{\mathbf v}|^2f\, \tn{d}{\xb} \tn{d}{\mathbf v}}_{ \text{ion energy}}
 + \underbrace{ k_\tB T_\te \int n_\te \ln n_\te\, \tn{d}{\xb}}_{\text{electron energy}}
 + \underbrace{\frac{1}{2\mu_0} \int |{\Bb}|^2\, \tn{d}{\xb}}_{\text{magnetic energy}}.
\ee 
We remark that quasi-neutrality $n = n_\te$ makes the electron continuity equation in~\eqref{model} redundant, since it is guaranteed by taking the velocity moment for Vlasov equation automatically. And keeping this equation brings the needs to conserve the quasi-neutrality relation and the positivity of $n_e$ numerically.

Several numerical methods and codes have already been proposed for this class of models, such as the current advance method~\cite{current2D}, based on which a Particle-In-Cell code CAMELIA~\cite{CAMELIA} and an Eulerian code~\cite{valentini} are developed; Pegasus~\cite{Pegasus}, in which a constrained transport method is used to enforce the divergence-free constraint on the magnetic field; and CHIEF~\cite{CHIEF}, in which an inertial electron fluid equation without approximation is adopted. Also in~\cite{densityjapan}, the low density issue is overcome by modifying the model to include finite electron inertia effects. The reader can find a more complete review about hybrid simulations in~\cite{YinL}.
Because the energy~\eqref{H} is non-quadratic, these methods do not provide energy conservation, which is important for long time simulations, for example in understanding the role of parametric instabilities in turbulence generation and proton heating in~\cite{para}. Related works also include \cite{chacon1, chacon2}, in which, for a modified kinetic-ion fluid-electron hybrid model with a quadratic energy, energy and momentum conserving schemes (on curvilinear meshes) are constructed. However, the conservation of the relation $p_e = T_e n_e^\gamma$ about the density and pressure of electrons is not addressed. 
Structure-preserving discretizations of this hybrid model represented with vector potentials are investigated in~\cite{potential}, in which the distribution functions depend on the canonical momentum (not velocity). 

Our discretizations follow the recent developments of structure-preserving methods~\cite{Feng, HLW}, which have been proposed with the purpose of preserving the geometric properties (such as symplectic structure) of the given systems and have very good long time behaviors. In plasma physics, some structure-preserving methods~\cite{Qin2, Qin3, Qin4, Qin5, GEMPIC, curviSISC, Morrison, 3, Li} have been proposed for kinetic equations and hybrid models. Specifically, in this work, following~\cite{3, GEMPIC}, we regard magnetic field as a 1-form, and corresponding finite element spaces are specially chosen (with B-splines as basis functions) such that  the commuting diagram proposed in~\cite{FEEC} holds. There are several big differences between our work and~\cite{GEMPIC, 3}. The main difference is that there are higher nonlinearities in our model. For example, in the anti-symmetric bracket we introduce in the following, we have integrals of products of three to six functions, for which numerical quadratures have to be used. Moreover, density is represented as the sum of finite number of particles, so smoothed delta functions are necessary to use.

There is a  Poisson bracket proposed in~\cite{Tronci} when $n_e$ is regarded as an independent unknown for the hybrid model. In this work, we propose an anti-symmetric bracket and numerical methods for the hybrid model with unknowns ($f, {\mathbf B}$), i.e., $n_e$ is not an independent unknown (it is replaced by $q/e\int f \mathrm{d}{\mathbf v}$) and $\frac{\partial n_\te}{\partial t} +  \nabla \cdot \left(n_\te \mathbf u  \right) = 0$ is not solved numerically. 
Time discretizations are derived by using Poisson splitting, i.e., by splitting the anti-symmetric bracket into several anti-symmetric sub-brackets, and get corresponding sub-steps. For the sub-step related with the non-quadratic term (thermal energy of electrons) in the Hamiltonian, we use the Hamiltonian splitting methods (for cartesian coordinate case) and discrete gradient methods (for curvilinear coordinate case) to conserve the energy in high precision or exactly.  
Numerical methods are constructed and detailed for the isothermal electron case, which can also be applied to the  adiabatic electron case. Note that as the anti-symmetric bracket is not Poisson, there is no Poisson structure preserved by the discretizations. However, total energy, divergence free property of magnetic field and positivity of density are conserved by the discretizations based on this bracket.

The paper is organized as follows. In section~\ref{sec:2}, an anti-symmetric bracket is introduced for the formulation with unknowns $(f, {\mathbf B})$, and time discretizations are done using Poisson splitting based on this anti-symmetric bracket. In section~\ref{sec:4}, curvilinear coordinates and the bracket in terms of differential form unknowns,  finite element methods and Particle-In-Cell methods are introduced, based on which discrete Poisson brackets are obtained by discrete functional derivatives. In section~\ref{sec:6}, computational details for sub-steps are discussed, also conservations of energy and divergence free property of magnetic field are discussed in details. In section~\ref{sec:7} and~\ref{sec:8}, performances and validations are done for the numerical schemes constructed. In section~\ref{sec:conclusion}, we conclude the paper with a summary and an outlook to future works.

\section{Poisson brackets}\label{sec:2}
Let us first get rid of physical constants by introducing the following units:
\be \label{units}
\begin{gathered}
t = \frac{t'}{\Omega}\,, \qquad \mathbf v = \mathbf v' v_\tA\,,
\qquad \xb = \xb' \frac{v_\tA}{\Omega}\,,\qquad \Bb = \Bb' B_0\,,\qquad \mathbf E = \mathbf E' v_A B_0\,,
\\[2mm]
n_e = n'_e n_0\,,\qquad \mathbf u = \mathbf u' v_\tA\,,\qquad \kappa := \frac{k_\tB T_\te}{m v_\tA^2}\,,\qquad \mathbf J = \mathbf J' e n_0 v_\tA\,,\qquad \mathcal H = \mathcal H' m n_0 v_\tA^2\,.
\end{gathered}
\ee
Here, primed quantities are dimension-less, $\Omega = qB_0/m$ denotes the ion cyclotron frequency, $B_0$ is the characteristic magnetic field strength and $v_\tA = B_0/\sqrt{\mu_0 m n_0}$ stands for the Alfv\'en velocity with $n_0$ the characteristic particle density. The parameter $\kappa$ is the ratio of electron energy to Alfv\'en kinetic energy. After normalization, we get the following dimensionless equations for the case $q = e$,
\begin{equation}
\label{eq:dimensionlessmodel}
\begin{aligned}
& \frac{\partial f}{\partial t} + {\mathbf v} \cdot \frac{\partial f}{\partial \xb} +({\mathbf E} + {\mathbf v} \times {\Bb}) \cdot \frac{\partial f}{\partial {\mathbf v}} = 0\,,\\
&\frac{\partial \Bb}{\partial t} = - \nabla \times {\mathbf E}\,,\quad \nabla \cdot \Bb = 0\,,\\
& \frac{\partial n_\te}{\partial t} +  \nabla \cdot \left(n_\te \mathbf u  \right) = 0, \quad {\mathbf u} = {\int {\mathbf v}f \mathrm{d}{\mathbf v}} \Big/ {\int f \mathrm{d}{\mathbf v} }, \\
& {\mathbf E} = - \kappa \frac{\nabla n_\te}{n_\te} - \left(\mathbf u - \frac{\mathbf J}{n_\te}  \right) \times {\Bb}, \quad \mathbf J = \nabla \times {\mathbf B}.\\
\end{aligned}
\end{equation}
It is interesting to note that the model \eqref{model}-\eqref{H}, being a combination of well-known plasma equations, is a non-canonical Hamiltonian system, i.e., it can be written in the formulation of Poisson brackets proposed in~\cite{Tronci}, in which $n_e$ is regarded as an independent unknown. Our aim is to work directly on the equations with only two unknowns $(f, {\mathbf B})$ and do not solve $\frac{\partial n_\te}{\partial t} +  \nabla \cdot \left(n_\te \mathbf u  \right) = 0$ numerically. For this, we propose an anti-symmetric bracket $\{\cdot,\cdot\}: V\times V \to V$, where $V$ denotes the vector space of functionals $\mathcal{P}: (f, \Bb) \mapsto \mathbb R$ on the space of unknowns, of the form 
\be \label{bra0}
\{ \mathcal{P},  \mathcal{Q} \} =  \{  \mathcal{P},  \mathcal{Q} \}_{vv}
 + \{  \mathcal{P},  \mathcal{Q} \}_{Bv} 
 + \{  \mathcal{P},  \mathcal{Q} \}_{BB}+ \{  \mathcal{P},  \mathcal{Q} \}_{xv}\,,
\ee
where\footnote{We shall omit the primes on normalized quantities in what follows and write everything in the units \eqref{units}.}
\begin{subequations} \label{bras}
\begin{align}
\{  \mathcal{P},  \mathcal{Q} \}_{vv} &=  -\int  \frac{1}{n}  f({\mathbf v}) f({\mathbf v}')  {\Bb} \cdot \left( \nabla_v \frac{\delta  \mathcal{P}}{\delta f} \times \nabla_{v'} \frac{\delta \mathcal{Q}}{\delta f} \right)\,  \tn{d}{\xb}  \tn{d}{\mathbf v} \tn{d}{\mathbf v}' + \int f  {\Bb} \cdot \nabla_v \frac{\delta  \mathcal{P}}{\delta f} \times \nabla_v \frac{\delta \mathcal{Q}}{\delta f}\,  \tn{d}{\xb}  \tn{d}{\mathbf v}\,,  \label{bras:3}
 \\[1mm]
 \{  \mathcal{P},  \mathcal{Q} \}_{Bv} &= \int \frac{f}{n}  \mathbf{B} \cdot \left[ \nabla_v \frac{\delta  \mathcal{P}}{\delta f} \times \left(\nabla \times \frac{\delta \mathcal{Q}}{\delta {\Bb}}\right) -   \nabla_v \frac{\delta \mathcal{Q}}{\delta f} \times \left(\nabla \times \frac{\delta  \mathcal{P}}{\delta {\Bb}}\right)\right] \tn{d}{\xb}  \tn{d}{\mathbf v}\,,  \label{bras:4}
 \\[1mm]
 \{  \mathcal{P},  \mathcal{Q} \}_{BB} &= -\int \frac{1}{n}  {\Bb} \cdot \left(\nabla \times \frac{\delta  \mathcal{P}}{\delta {\Bb}}\right) \times \left(\nabla \times \frac{\delta \mathcal{Q}}{\delta {\Bb}}\right)  \tn{d}{\xb}\,,  \label{bras:5}
 \\[1mm]
 \{  \mathcal{P},  \mathcal{Q} \}_{xv} &= \int f  \left(\nabla_x\frac{\delta  \mathcal{P}}{\delta f} \cdot \nabla_v \frac{\delta \mathcal{Q}}{\delta f}
 - \nabla_v\frac{\delta  \mathcal{P}}{\delta f} \cdot \nabla_x \frac{\delta \mathcal{Q}}{\delta f}\right)  \tn{d}{\xb}  \tn{d}{\mathbf v} \,.  \label{bras:1}
\end{align}
\end{subequations}
Note that $n = \int f {\mathrm d}{\mathbf v}$ and $n_e$ is absent from this bracket. 
The Hamiltonian (total energy) can be written as 
\be \label{H2}
 \cH =  \frac{1}{2} \int |{\mathbf v}|^2f\, \tn{d}{\xb} \tn{d}{\mathbf v}
 + \frac{1}{2} \int |{\Bb}|^2\, \tn{d}{\xb}\, +  \kappa \int \left( \int f\,\tn d \vb\right) \ln \left( \int f\,\tn d \vb\right) \tn{d}{\xb},
\ee
which is the sum of the ion kinetic energy, the magnetic Hamiltonian, and electron thermal energy. The functional derivatives about $f$ and ${\mathbf B}$ read
$$
\del{\cH}{f} = \frac{|\vb|^2}{2} + \kappa \left(1 + \ln \left( \int f\,\tn d \vb\right)\right) \,,\qquad \del{\cH}{\Bb} = \Bb \,.
$$
Based on the above anti-symmetric bracket and Hamiltonian~\eqref{H2}, we can reformulate the hybrid model as 
\begin{equation}\label{eq:zform}
\dot{\mathcal{Z}} = \{ \mathcal{Z}, \mathcal{H}  \}, \quad \mathcal{Z} = (f, {\mathbf B}).
\end{equation}
The advantages of this formulation and bracket~\eqref{bras} include that 1) the
quasi-neutrality relation of ions and electrons is built in, as $n_e$ is not an independent unknown; 2) the positivity of $n$ is guaranteed numerically, as it is written as the sum of finite number of particles after Particle-In-Cell discretizations; 3) the anti-symmetric bracket~\eqref{bras} and formulation~\eqref{eq:zform} build the foundation to use discrete gradient methods~\cite{Gonzalez, DIS} for conserving energy numerically,  although the anti-symmetric bracket~\eqref{bras} does not satisfy Jacobi identity. 

%\section{Splitting schemes}\label{sec:3}

% Using the parametrized functionals 
% \begin{align}
%  \cF_f(t,\xb',\vb') &= \int f(t, \xb,\vb)\, \delta^3(\xb - \xb')\delta^3(\vb - \vb')\, \tn d \xb \tn d \vb
%  \\[1mm]
%  \cF_B(t,\xb') &= \int \Bb(t,\xb)\, \delta^3(\xb - \xb')\,\tn d \xb\,,  \label{FB:delta}
% \end{align} 
% the model equations are obtained from
% $$
%  \frac{\partial \cF_i}{\partial t} = \{\cF_i, \cH\} \,, \quad i = f, {\mathbf B}.
% $$
In the numerical discretizations we will use Poisson splitting methods \cite{Nicolas, Qin4, 2020en} based on the decomposition 
of the full bracket \eqref{bra0} into the four parts \eqref{bras}, which leads to the following 4 sub-steps conserving energy, mass $\int f \mathrm{d}{\mathbf v}\mathrm{d}{\mathbf x}$, and momentum $\int {\mathbf v} f \mathrm{d}{\mathbf v}\mathrm{d}{\mathbf x}$:
\begin{alignat}{2}
      \tn{sub-step}\, vv: \qquad & \frac{\partial  \mathcal{P}}{\partial t} = \{ \mathcal{P}, \cH\}_{vv}  \quad \implies\quad 
    \pder{f}{t} = -(\vb \times \Bb + \Bb \times \frac{\int \vb f\,\tn d\vb}{n}) \cdot \nabla_v f\,,  \label{stepvv}
    \\[2mm]
    \tn{sub-step}\ bv: \qquad & \frac{\partial  \mathcal{P}}{\partial t} = \{\mathcal{P}, \cH\}_{Bv} \quad \implies\quad 
    \left\{
    \begin{aligned}
    &\pder{f}{t} = -\frac{1}{n} (\nabla \times \Bb) \times \Bb \cdot \nabla_v f\,,
    \\
    &\pder{\Bb}{t} = -\nabla \times \left( \Bb \times \frac{1}{n} \int \vb f\,\tn d\vb \right)\,,
    \end{aligned}
    \right.  \label{stepbv}
    \\[2mm]
    \tn{sub-step}\ bb: \qquad & \frac{\partial \mathcal{P}}{\partial t} = \{\mathcal{P}, \cH\}_{BB} \quad \implies\quad 
    \pder{ {\mathbf B}}{t} = - \nabla \times \left( \frac 1n \nabla \times {\mathbf B} \times {\mathbf B}\right)\,,  \label{stepbb}
    \\[2mm]
    \tn{sub-step}\ xv: \qquad & 
    \begin{aligned}
    &\frac{\partial \mathcal{P}}{\partial t} = \{\mathcal{P}, \cH\}_{xv}
     \quad \implies 
    \pder{f}{t} = - \mathbf v \cdot \nabla_x f + \kappa \frac{\nabla n}{n} \cdot \nabla_v f\,, 
    \end{aligned}  \label{sys:formIII}
\end{alignat}
where $n = {\int  f\,\tn d \vb}$.

\section{Discretization of brackets}\label{sec:4}
In this section, the anti-symmetric bracket~\eqref{bras} for the hybrid model is discretized following three steps. Firstly, the bracket is presented in terms of curvilinear coordinates by means of differential forms; then finite element methods~\cite{FEEC, GEMPIC, buffa} and Particle-In-Cell methods are used for the discretizations; finally, functional derivatives are discretized to give the discrete bracket. 

\subsection{Curvilinear coordinates and differential forms}

Here we write the considered model in curved space in terms of differential forms \cite{Frankel}. We use the same notations as in \cite{3, polar}.
The smooth, invertible coordinate transformation (or "mapping") is ${\Fb}: \hat{\Omega} \rightarrow {\Omega}, {\etab} \rightarrow {\Fb}({\etab}) = {\xb}$, from the logical cuboid $\hat{\Omega} = [0,1]^3$ to the physical domain $\Omega \subset \mathbb{R}^3$. Moreover, ${\etab} = (\eta_1, \eta_2, \eta_3) \in \hat{\Omega}$ and ${\xb} = (x,y,z) \in \Omega$ are the logical and cartesian coordinates, respectively.  The Jacobian matrix of this coordinate transformation is
$$
DF: \hat{\Omega} \rightarrow \mathbb{R}^{3 \times 3}, \quad (DF)_{i,j} = \frac{\partial F_i}{\partial \eta_j}.
$$
For convenience, we also introduce the following notations,
$$
\hat{\nabla} = (\partial_{\eta_1}, \partial_{\eta_2},  \partial_{\eta_3})^\top\,,\quad  {\nabla} = (\partial_{x}, \partial_{y},  \partial_{z})^\top\,,\quad  G = DF^\top DF\,,\quad \ g = \text{det}G = \text{det}(DF)^2.
$$
We shall frequently use the identity 
$
 M\ab \times M\bb = \det (M) M^{-\top} (\ab \times \bb)
$
for matrix $M$ of size $3 \times 3$ to transform cross products. 
Scalar fields transform as $a({\Fb}({\etab} )) = \hat{a}({\etab} )$, the components of a vector-field ${\mathbf a} \in T_{\xb}\Omega$  in curvilinear (logical) coordinates, denoted by $\hat{\boldsymbol a}$, are defined by the relation ${\boldsymbol a}({\Fb}({\etab})) = DF({\etab})\hat{\boldsymbol a}({\etab})$. The relations of gradient, curl, and divergence in logical and physical coordiantes are: 
\begin{equation} \label{transf:d}
\nabla a = DF^{-\top} \hat{\nabla}\hat{a}, \quad \nabla \times {\boldsymbol a} = \frac{1}{\sqrt{g}}DF \hat{\nabla} \times (G\hat{{\boldsymbol a}}), \quad \nabla \cdot {\boldsymbol a} = \frac{1}{\sqrt{g}}\hat{\nabla}\cdot (\sqrt{g}\hat{\boldsymbol a}).
\end{equation}
Scalar fields $\hat{a} = \hat{a}({\etab})$ and vector fields $a \in T\Omega$ with components $\hat{\boldsymbol a} = \hat{\boldsymbol a}({\etab})$ can be related to differential $p$-forms $a^p \in \Lambda^p(\Omega): T\Omega \times \cdots \times T\Omega \rightarrow \mathbb{R}, p \in \{0,1,2,3\}$ with components $\hat{ a}^p = \hat{ a}^p({\etab} )$ in the following way:
\begin{equation} \label{proxies}
\begin{aligned}
& a^0 = \hat{a}^0  \leftrightarrow \hat{a}^0 = \hat{a},\\
& a^1 = \hat{a}_1^1 \tn{d}\eta^1 + \hat{a}_2^1 \tn{d}\eta^2 + \hat{a}_3^1 \tn{d}\eta^3 \leftrightarrow \hat{\boldsymbol a}^1 = \left( \begin{matrix}  \hat{a}_1^1 \\ \hat{a}_2^1 \\ \hat{a}_3^1
\end{matrix} \right) = G \hat{\boldsymbol a},\\
& a^2 = \hat{a}_1^2 (\tn{d}\eta^2 \wedge \tn{d}\eta^3) + \hat{a}_2^2 (\tn{d}\eta^3 \wedge \tn{d}\eta^1)  + \hat{a}_3^2 \tn{d}\eta^1 \wedge \tn{d}\eta^2 \leftrightarrow  \hat{\boldsymbol a}^2 = \left( \begin{matrix}  \hat{a}_1^2 \\ \hat{a}_2^2 \\ \hat{a}_3^2
\end{matrix} \right) = \sqrt{g} \hat{\boldsymbol a},\\
& a^3 = \hat{a}^3 (\tn{d}\eta^1 \wedge \tn{d}\eta^2 \wedge \tn{d}\eta^3)  \leftrightarrow \hat{a}^3 = \sqrt{g} \hat{a}.
\end{aligned}
\end{equation}
The functions $\hat a^p : \hat\Omega \to \mathbb R$ can be thought of as proxies of the differential forms $a^p$. For functionals acting on the space of proxies we use the push-forward of $p$-forms to define 
\begin{align*}
    & \mathcal{P}(a(\xb)) = \mathcal{P}(\hat a^0(\etab)) &&=: \hat{\mathcal{P}}(\hat a^0(\etab))\,,
    \\[2mm]
    &\mathcal{P}({\mathbf a}({\xb})) = \mathcal{P}( DF^{-\top} \hat{\mathbf a}^1({\etab})) &&=: \hat{\mathcal{P}}(\hat{\mathbf a}^1({\etab}))\,,
    \\[2mm]
    &\mathcal{P}({\mathbf a}({\xb})) = \mathcal{P}( DF/\sqrt g\, \hat{\mathbf a}^2({\etab})) &&=: \hat{\mathcal{P}}(\hat{\mathbf a}^2({\etab}))\,,
    \\[2mm]
    & \mathcal{P}(a(\xb)) = \mathcal{P}(\hat a^3(\etab)/\sqrt g) &&=: \hat{\mathcal{P}}(\hat a^3(\etab))\,.
\end{align*}
The corresponding functional derivatives can be computed by duality using again the push-forward (here for a 1-form)
\be \label{dF}
     \delta \mathcal{P}({\ab}({\xb})) = \int_{\Omega} \frac{\delta \mathcal{P}}{\delta {\ab}} \cdot \delta{\ab} \,\tn{d}{\xb} = \int_{\hat\Omega} \frac{\delta \mathcal{P}}{\delta {\ab}}({\Fb}({\etab})) \cdot DF^{-\top} \delta \hat{\ab}^1 \sqrt{g} \,\tn{d}{\etab} = \int_{\hat{\Omega}} \frac{\delta \hat{\mathcal{P}}}{\delta {\hat{\ab}}^1} \cdot \delta{\hat{\ab}}^1 \tn{d}{\etab}\,,
\ee
and similar for the other $p$-forms, which leads to
\be\label{eq:rreela}
\frac{\delta \mathcal{P}}{\delta {a}} = \frac{1}{\sqrt{g}}  \frac{\delta \hat{\mathcal{P}}}{\delta { \hat{a}^0}}\,, 
\qquad
\frac{\delta \mathcal{P}}{\delta {\ab}} = \frac{1}{\sqrt{g}}  DF \frac{\delta \hat{\mathcal{P}}}{\delta \hat{\ab}^1}\,,
\qquad
\frac{\delta \mathcal{P}}{\delta {\ab}} = DF^{-\top} \frac{\delta \hat{\mathcal{P}}}{\delta \hat{\ab}^2}\,,
\qquad 
\frac{\delta \mathcal{P}}{\delta {a}} = \frac{\delta \hat{\mathcal{P}}}{\delta { \hat{a}^3}}\,. 
\ee
Hence the functional derivative with respect to a $p$-form (proxy) transforms as a $3$-$p$-form (proxy). This is necessary for integrals as in \eqref{dF} to be well-defined, and stems from the fact that functional derivatives live in the dual spaces.

The brackets \eqref{bras} can be written in terms of differential forms by applying the transformations \eqref{transf:d}, \eqref{proxies} and \eqref{eq:rreela}. It is thus necessary to make the choice of $p$-form representation for the unknowns $(f,\Bb)$. The choice is free in principle but can be informed by boundary conditions or by certain properties of the ensuing numerical schemes. There are two different choices which lead to a weak and a strong form of the field equations, respectively. This translates to weak or strong conservation laws. The "weak scheme" has the magnetic field $\Bb$ as a 1-form; the "strong scheme" has the magnetic field $\Bb$ as a 2-form. The distribution function $f$ is a 3-form in both cases. In this work, we choose a weak formulation to construct numerical methods, the strong formulation can be considered similarly. Specifically, the sub-brackets in terms of differential forms are
\begin{subequations}
\small{
\begin{align*}
\{  \hat{\mathcal{P}},  \hat{\mathcal{Q}} \}_{vv} &= \int  \hat f^3  DF^{-\top}\hat\Bb^1 \cdot \nabla_v \frac{\delta  \hat{\mathcal{P}}}{\delta \hat f^3} \times \nabla_v \frac{\delta  \hat{\mathcal{Q}}}{\delta \hat f^3} \, \tn{d}{\etab}  \tn{d}{\mathbf v}\\
& -\int  \frac{1}{\hat{n}^3}  \hat f^3({\mathbf v}) \hat f^3({\mathbf v}')  DF^{-\top}\hat{\Bb}^1 \cdot \left( \nabla_v \frac{\delta  \hat{\mathcal{P}}}{\delta \hat f^3}  \times  \nabla_{v'} \frac{\delta  \hat{\mathcal{Q}}}{\delta \hat f^3}  \right)\,  \tn{d}{\etab} \tn{d}{\mathbf v} \tn{d}{\mathbf v}'\,, 
 \\[1mm]
 \{  \hat{\mathcal{P}},  \hat{\mathcal{Q}} \}_{Bv} &= \int \frac{\hat f^3}{\hat{n}^3} DF^{-\top}\hat \Bb^1 \cdot \left[ \nabla_v \frac{\delta  \hat{\mathcal{P}}}{\delta \hat f^3} \times \left(DF\hat\nabla \times \frac{G}{\sqrt g} \frac{\delta  \hat{\mathcal{Q}}}{\delta \hat{\Bb}^1} \right) -   \nabla_v \frac{\delta  \hat{\mathcal{Q}}}{\delta \hat f^3} \times \left(DF\hat\nabla \times \frac{G}{\sqrt g}\frac{\delta  \hat{\mathcal{P}}}{\delta \hat{\Bb}^1}\right)\right] \, \tn{d}{\etab}  \tn{d}{\mathbf v}\,, 
 \\[1mm]
 \{ \hat{\mathcal{P}},  \hat{\mathcal{Q}} \}_{BB} &= -\int \frac{1}{\hat{n}^3}  DF^{-\top}\hat{\Bb}^1 \cdot \left(DF \hat \nabla \times \frac{G}{\sqrt g} \frac{\delta  \hat{\mathcal{P}}}{\delta \hat{\Bb}^1}\right) \times \left(DF \hat\nabla \times \frac{G}{\sqrt g} \frac{\delta  \hat{\mathcal{Q}}}{\delta \hat{\Bb}^1}\right)\, \tn{d}{\etab}\,, 
 \\[1mm]
 \{  \hat{\mathcal{P}},  \hat{\mathcal{Q}} \}_{xv} &= \int \hat f^3  \left( DF^{-\top} \hat\nabla_\eta \frac{\delta  \hat{\mathcal{P}}}{\delta \hat f^3}  \cdot \nabla_v \frac{\delta  \hat{\mathcal{Q}}}{\delta \hat f^3} 
 - \nabla_v\frac{\delta  \hat{\mathcal{P}}}{\delta \hat f^3} \cdot DF^{-\top}\hat \nabla_\eta  \frac{\delta  \hat{\mathcal{Q}}}{\delta \hat f^3} \right) \tn{d}{\etab}  \tn{d}{\mathbf v} \,,
\end{align*}
}
\end{subequations}
where $\hat{n}^3 = \int \hat{f}^3 \mathrm{d}{\mathbf v}$.

\subsection{Phase-space discretization}\label{sec:5}
As for the discretizations in space (phase-space), we choose the finite element method for the unknowns of fields, and the Particle-In-Cell (PIC) method for the distribution function. \\

\noindent{\bf Commuting diagram with B-splines.}
We perform the spatial discretizations in the framework of Finite Element Exterior Calculus (FEEC)~\cite{FEEC, GEMPIC, buffa}.
Finite element (FE) spaces and corresponding projectors are chosen such that the following diagram commutes,
\begin{equation}
\large
\begin{aligned}
\label{diagram}
 \xymatrix{
    H^{1}(\hat{\Omega}) \ar[rr]^{\hat{\nabla}} \ar[d]_{\Pi_0} && H(\text{curl}, \hat{\Omega}) \ar[rr]^{\hat{\nabla}\times} \ar[d]^{\Pi_1} && H(\text{div}, \hat{\Omega}) \ar[rr]^{\hat{\nabla}\cdot} \ar[d]^{\Pi_2}  && L^2(\hat{\Omega}) \ar[d]^{\Pi_3} \\
      V_0 \ar[rr]^{\hat{\nabla}}  && V_1 \ar[rr]^{\hat{\nabla}\times}  && V_2  \ar[rr]^{\hat{\nabla}\cdot} && V_3 
}
\end{aligned}
\end{equation}
where $V_0, V_1, V_2$ and $V_3$ are finite element spaces in which fields (proxies of $p$-forms, $p\leq 3$) are discretized. The projectors $\Pi_n, n = 0, 1, 2, 3$ are based on inter-/histopolation at/between Greville points of the B-splines which span the finite element spaces. For more details we refer to \cite{3} which uses exactly the same basis functions and projectors. The finite element spaces are written as
\begin{alignat*}{2}
& V_0 := \text{span} \{\Lambda^0_i | 0 \le i < N_0  \}, &&\ni \hat{a}^0_h(t, {\etab}) = \sum_{i=0}^{N_0-1} a_i(t)\Lambda^0_i({\etab}) \\
& V_1 := \text{span}  \left\{
\left(\begin{matrix}
  \Lambda^1_{1,i} \\
    0\\
      0
\end{matrix} \right),
\left(\begin{matrix}
  0\\
     \Lambda^1_{2,i} \\
      0
\end{matrix} \right),
\left(\begin{matrix}
  0 \\
    0\\
       \Lambda^1_{3,i} 
\end{matrix} \right)
\Bigg| 
\begin{matrix}
  0 \le i < N^{1}_1\\
    0 \le i < N^{1}_2 \\
      0 \le i < N^{1}_3
\end{matrix} 
\right\}
&&\ni \hat{\mathbf a}^1_h(t, {\etab}) = \sum_{\mu = 1}^3 \sum_{i=0}^{N_\mu^1-1} a_{\mu,i}(t) \Lambda^1_{\mu, i}({\etab}) {\mathbf e}_\mu,\\
& V_2 := \text{span}  \left\{
\left(\begin{matrix}
  \Lambda^2_{1,i} \\
    0\\
      0
\end{matrix} \right),
\left(\begin{matrix}
  0\\
     \Lambda^2_{2,i} \\
      0
\end{matrix} \right),
\left(\begin{matrix}
  0 \\
    0\\
       \Lambda^2_{3,i} 
\end{matrix} \right)
\Bigg| 
\begin{matrix}
  0 \le i < N^{2}_1\\
    0 \le i < N^{2}_2 \\
      0 \le i < N^{2}_3
\end{matrix} 
\right\}
&&\ni \hat{\mathbf a}^2_h(t, {\etab}) = \sum_{\mu= 1}^3 \sum_{i=0}^{N_\mu^1-1} a_{\mu,i}(t) \Lambda^2_{\mu, i}({\etab}) {\mathbf e}_\mu,\\
& V_3 := \text{span} \{\Lambda^3_i | 0 \le i < N_3  \},  &&\ni\hat{a}^3_h(t, {\etab}) = \sum_{i=0}^{N_3-1} a_i(t)\Lambda^3_i({\etab}) \,.
\end{alignat*}
Here, the functions $\Lambda_i^n : \hat\Omega \to \mathbb R$ are tensor products of uni-variate B-splines of different degrees, as described in \cite{3, buffa, GEMPIC}, and ${\mathbf e}_1 = (1,0,0)^\top$, ${\mathbf e}_2 = (0,1,0)^\top$, ${\mathbf e}_3 = (0,0,1)^\top$. The dimensions are
$$
 \tn{dim} V_0 = N_0\,,\quad \tn{dim} V_1 = N_1 =  \sum_{d=1}^3 N^1_d\,, \quad \tn{dim} V_2 = N_2 =  \sum_{d=1}^3 N^2_d\,,\quad \tn{dim} V_3 = N_3\,.
$$
To simplify the notation, the finite element coefficients $a_i$ and basis functions are stacked in column vectors, e.g. ${\mathbf a} := (a_i)_{0 \le i  <N_0} \in \mathbb{R}^{N_0}$, and ${{\mathbf \Lambda}^0} := (\Lambda^0_i)_{0\le  i  < N_0} \in \mathbb{R}^{N_0}$. Spline functions can then be compactly written as  
\begin{equation*}
\begin{aligned}
& \hat{a}_h^0 = {\mathbf a}^{\top} {\mathbf  \Lambda}^0,\\
& (\hat{\mathbf a}_h^1)^\top = (\underbrace{a_{1,0}, \cdots, a_{1,N^1_1-1}}_{=: {\mathbf a}_1^\top}, \underbrace{a_{2,0}, \cdots, a_{2,N^1_2-1}}_{=: {\mathbf a}_2^\top}, \underbrace{a_{3,0}, \cdots, a_{3,N^1_3-1}}_{=: {\mathbf a}_3^\top}) 
\left(\begin{matrix}
  {\mathbf \Lambda}^1_{1} & 0 & 0 \\
    0 & {\mathbf \Lambda}^1_{2}  & 0 \\
      0 & 0 & {\mathbf \Lambda}^1_{3}
\end{matrix} \right) =: {\mathbf a}^\top \mathbb{\Lambda}^1,\\
& (\hat{\mathbf a}_h^2)^\top = (\underbrace{a_{1,0}, \cdots, a_{1,N^2_1-1}}_{=: {\mathbf a}_1^\top}, \underbrace{a_{2,0}, \cdots, a_{2,N^2_2-1}}_{=: {\mathbf a}_2^\top}, \underbrace{a_{3,0}, \cdots, a_{3,N^2_3-1}}_{=: {\mathbf a}_3^\top}) 
\left(\begin{matrix}
  {\mathbf \Lambda}^1_{1} & 0 & 0 \\
    0 & {\mathbf \Lambda}^1_{2}  & 0 \\
      0 & 0 & {\mathbf \Lambda}^1_{3}
\end{matrix} \right) =: {\mathbf a}^\top \mathbb{\Lambda}^2,\\
& \hat{a}_h^3 = {\mathbf a}^{\top} {\mathbf  \Lambda}^3,
\end{aligned}
\end{equation*}
where $\mathbb{\Lambda}^1 \in \mathbb{R}^{N_1 \times 3}$ and $\mathbb{\Lambda}^2 \in \mathbb{R}^{N_2 \times 3}$.
In this setting the discrete representations of the exterior derivatives can be written as matrices acting on finite element coefficients, 
$$
V_1\ni \hat{\nabla} \hat{a}^0_{h} = ({\mathbb{G} {\mathbf {\mathbf a}}})^\top \mathbb{\Lambda}^1, \quad V_2 \ni  \hat{\nabla} \times \hat{\mathbf a}_h^1 = ({\mathbb{C} {\mathbf a}})^\top \mathbb{\Lambda}^2\,,\quad V_3 \ni \hat \nabla \cdot \hat \ab^2_h = (\mathbb D \ab)^\top \mathbb \Lambda^3\,.
$$
where $\mathbb{G} \in \mathbb{R}^{N_1 \times N_0}$,  $\mathbb{C} \in \mathbb{R}^{N_2 \times N_1}$ and $\mathbb{D} \in \mathbb{R}^{N_3 \times N_2}$ are sparse and contain only zeros and ones.
Finally, the (symmetric) mass matrices corresponding to the discrete spaces $V_n, n = 0, 1,2,3$ follow from the $L^2$-inner products and contain metric coefficients, according to the transformation rules of proxies in \eqref{proxies}~(push-forward):
\begin{align*}
&\mathbb{M}_0 := \int {\mathbf \Lambda^0} ({\mathbf \Lambda^0})^\top \sqrt{g} \,\tn{d}^3{\etab} &&\in \mathbb R^{N_0 \times N_0}\,,
\\
&\mathbb{M}_1 := \int {\mathbb{\Lambda}^1} G^{-1} ({\mathbb{\Lambda}^1})^\top \sqrt{g} \,\tn{d}^3{\etab} &&\in \mathbb R^{N_1 \times N_1}\,,
\\
&\mathbb{M}_2 := \int {\mathbb{\Lambda}^2} G ({\mathbb{\Lambda}^2})^\top \frac{1}{\sqrt{g}}\, \tn{d}^3{\etab} &&\in \mathbb R^{N_2 \times N_2}\,,
\\
&\mathbb{M}_3 := \int {\mathbf \Lambda^3} ({\mathbf \Lambda^3})^\top \frac{1}{\sqrt{g}}\, \tn{d}^3{\etab} &&\in \mathbb R^{N_3 \times N_3}\,.
\end{align*}
These mass matrices are sparse because of the compact supports of B-splines.

In this work, we choose a "weak" formulation of the field equations, which leads to the following discretizations:
\begin{align*}
\hat \Bb^1 \approx \hat \Bb^1_h = \bb^\top \mathbb \Lambda^1,
\end{align*}
i.e., the magnetic field is a 1-form. The functional derivatives are expressed in the range of  hodge-star operator on the spaces of unknowns. For example, the functional derivative about magnetic field (in  $V_1$) $\frac{\delta \hat{\mathcal{P}}}{\delta {\hat{\Bb}}^1}$ is discretized in $*V_1 := \{ \sqrt{g}G^{-1} {\mathbf a}^1 |  {\mathbf a}^1 \in  V_1 \}$. Then we have 
$$
 \delta{\mathbf B}_h^1 = \delta{\mathbf b} \cdot \mathbb \Lambda^1 \in V_1\,, \quad \frac{\delta \hat{\mathcal{P}}}{\delta {\hat{\Bb}}^1} = \partial\bb^\top \sqrt{g} G^{-1}\mathbb \Lambda^1 \in  *V_1,
$$
$$
 \delta \hat{\mathcal{P}}({\mathbf B}_h^1) =  \int_{\hat{\Omega}} \frac{\delta \hat{\mathcal{P}}}{\delta {\hat{\Bb}}^1} \cdot \delta{\hat{\Bb}}^1_h \tn{d}{\etab} = \partial\bb^\top \mathbb{M}_1 \delta{\mathbf b}.
$$
Viewing $\hat{\mathcal{P}}(\hat \Bb^1_h) =: P(\bb)$ as a functional on the discrete space $\mathbb R^{N_1}$, we have 
$$
 \delta \hat{\mathcal{P}}({\mathbf B}_h^1) =  \delta P({\mathbf b}) = \frac{\partial P}{\partial {\mathbf b}} \cdot \delta{\mathbf b},
$$
from which we get $\partial\bb = \mathbb{M}_1^{-1} \frac{\partial P}{\partial {\mathbf b}}$, and
\begin{equation} \label{eq:weakb}
\frac{\delta \hat{\mathcal{P}}}{\delta {\hat{\Bb}}^1} = (\sqrt{g} G^{-1}\mathbb \Lambda^1)^\top \mathbb{M}_1^{-1} \frac{\partial P}{\partial {\mathbf b}}\,.
\end{equation}

\noindent{\bf Particle-In-Cell method}
The kinetic distribution function is a 3-form in logical coordinates and is represented (discretized) by a finite number $K$ of  "markers", 
\begin{equation} \label{fh3}
\begin{aligned}
\hat{f}^3_h(t, {\etab}, {\mathbf v}) & = \sum_{k=1}^K w_k\, S({\etab} - {\etab}_k(t))\, \delta({\mathbf v} - {\mathbf v}_k(t)),
\end{aligned}
\end{equation}
Here, $\etab_k(t)$ and $\vb_k(t)$ are the position (in logical space $\hat\Omega$) and velocity of the $k$-th marker, also the solutions of the characteristics equations of the Vlasov equation. The number $w_k\in\mathbb R$ is the $k$-th marker's weight. We give the discretization of $\hat{n}^3$ as 
$$
\hat{n}^3_h = \sum_k w_k S(\etab - \etab_k(t)).
$$
The shape function $S:\hat\Omega \to \mathbb R$ can be chosen as the $\delta$-function~\cite{GEMPIC} or as a compactly supported function~\cite{martin}. The shape function is normalized, i.e. $\int S({\etab}) \tn{d}{\etab} = 1$. In this work the compactly supported smoothed shape functions are 
\begin{equation}\label{eq:shape_support}
S(\boldsymbol \eta) = \frac{1}{h_1h_2h_3} S_{k_1}\left(\frac{\eta_1}{h_1}\right)S_{k_2}\left(\frac{\eta_2}{h_2}\right)S_{k_3}\left(\frac{\eta_3}{h_3}\right), 
\end{equation}
where $S_k$ is defined as 
$$
S_0(\eta) := \mathbb{1}_{[-\frac{1}{2}, \frac{1}{2}]}, \quad S_k(\eta) = S_0 \star S_{k-1} = \int_{-\frac{1}{2}}^{\frac{1}{2}} S_{k-1}(\eta - y) \mathrm{d}{y}.
$$
Then we know that the localized support of $S({\boldsymbol \eta})$ is $\text{supp}(S) = [-h_1k_1, h_1k_1] \times  [-h_2k_2, h_2k_2] \times  [-h_3k_3, h_3k_3] $. Note that here $h_1, h_2, h_3$ may be chosen different from the cell size of fields' discretization. 

We can view the marker distribution function \eqref{fh3} as a parametrized functional of the trajectories $(\mathbf H(t), \mathbf V(t))$ (see the notation in~\eqref{notations}); its variation reads formally 
\be \label{deltaf3}
\begin{aligned}
 \delta \hat f^3_h(\mathbf H(t), \mathbf V(t)) = - \sum_{k=1}^K w_k \left[\delta {\boldsymbol  \eta}_k  \cdot  \hat{\nabla} S({\boldsymbol  \eta} - {\boldsymbol  \eta}_k)\, \delta({\mathbf v} - {\mathbf v}_k)
 + \delta {\mathbf  v}_k  \cdot S({\boldsymbol  \eta} - {\boldsymbol  \eta}_k) \nabla_v  \delta({\mathbf v} - {\mathbf v}_k) \right] \,.
\end{aligned}
\ee
Also we can view any functional $\hat{\mathcal{P}}$ of $\hat f_h^3$~\eqref{fh3} as a function  $P(\mathbf H(t), \mathbf V(t))$.  Then we have 
\begin{align}
\delta \hat{\mathcal{P}}(\hat f_h^3)&= \int \frac{\delta \hat{\mathcal{P}}}{\delta \hat{f}^3_h} \delta \hat f^3_h \, \tn{d}{\etab} \tn{d}{\mathbf v}\\
& =\sum_{k=1}^K w_k \left( \int S({\boldsymbol  \eta} - {\boldsymbol  \eta}_k) \hat\nabla \frac{\delta \hat{\mathcal{P}}}{\delta \hat{f}^3_h}\bigg|_{\vb_k} \cdot \delta {\boldsymbol  \eta}_k \tn{d}{\etab} + \int S({\boldsymbol  \eta} - {\boldsymbol  \eta}_k) \nabla_v \frac{\delta \hat{\mathcal{P}}}{\delta \hat{f}^3_h}\bigg|_{\vb_k} \cdot \delta {\vb}_k \tn{d}{\etab} \right)  \,.  \nonumber
\end{align}
On the other hand, for $P(\mathbf H(t), \mathbf V(t))$, we have
$$
 \delta P(\mathbf H(t), \mathbf V(t)) = \sum_{k=1}^K \left[ \frac{\partial P}{\partial {\etab}_k} \cdot \delta \etab_k + \frac{\partial P}{\partial \vb_k} \cdot \delta \vb_k \right].
$$
Therefore,
\begin{alignat}{2}
  \frac{\partial P}{\partial {\etab}_k} &= w_k\int S({\boldsymbol  \eta} - {\boldsymbol  \eta}_k) \hat\nabla \frac{\delta \hat{\mathcal{P}}}{\delta \hat{f}^3_h}\bigg|_{\vb_k} \tn{d}{\etab} && \approx w_k \hat\nabla \frac{\delta \hat{\mathcal{P}}}{\delta \hat{f}^3_h}\bigg|_{\etab_k,\vb_k}\,,  \label{fder:eta}
  \\[1mm]
  \frac{\partial P}{\partial {\mathbf v}_k} &= w_k\int S({\boldsymbol  \eta} - {\boldsymbol  \eta}_k) \nabla_v \frac{\delta \hat{\mathcal{P}}}{\delta \hat{f}^3_h}\bigg|_{\vb_k} \tn{d}{\etab} && \approx w_k \nabla_v \frac{\delta \hat{\mathcal{P}}}{\delta \hat{f}^3_h}\bigg|_{\etab_k,\vb_k} \,.   \label{fder:v}
\end{alignat}
The '$\approx$' in the above would become '$=$' when $S$ is a delta function as~\cite{GEMPIC}.

\noindent{\bf Notations}
 In order to write the discrete equations in the following sections in a more compact matrix-vector form, we introduce the following vectors and matrices:
\begin{equation}
\label{notations}
\begin{aligned}
&{\mathbf H} := (\eta_{1,1},  \cdots, \eta_{K,1}, \eta_{1,2}, \cdots, \eta_{K,2}, \eta_{1,3}, \cdots, \eta_{K,3})^\top \quad && \in \mathbb{R}^{3K}\,,
\\[2mm]
& {\mathbf V} :=   (v_{1,1},  \cdots, v_{K,1}, v_{1,2}, \cdots, v_{K,2}, v_{1,3}, \cdots, v_{K,3})^\top \quad && \in \mathbb{R}^{3K}\,,
\\[2mm]
&\mathbb{P}^n_\mu({\mathbf H}) := (\Lambda^n_{\mu,i}({\etab}_k))_{0\le i < N^n_\mu, 1 \le k \le K} \quad (n \in \{1, 2\}, \mu \in \{ 1, 2, 3\}) \quad && \in \mathbb{R}^{N_\mu^n \times K}\,,
\\[2mm]
& \mathbb{P}^n({\mathbf H}) := \text{diag}(\mathbb{P}^n_1, \mathbb{P}^n_2, \mathbb{P}^n_3), n \in \{1, 2\} \quad && \in \mathbb{R}^{N^n \times 3K}\,,
\\[2mm]
& \bar{G}^{-1}_{ab}({\mathbf H}) := \text{diag}(G^{-1}_{ab}({\etab}_1), \cdots, G^{-1}_{ab}({\etab}_{K})) \quad (a, b \in \{1,2,3 \}) \quad && \in \mathbb{R}^{K\times K}\,,
\\[2mm]
& \bar{G}^{-1} := ( \bar{G}^{-1}_{ab})_{1\le a, b \le 3} \quad && \in \mathbb{R}^{3K\times 3K}\,,
\\[2mm]
&  \bar{DF}^{-1}_{ab}({\mathbf H}) :=\text{diag}(DF^{-1}_{ab}({\etab}_1), \cdots, DF^{-1}_{ab}({\etab}_{K})) \quad (a, b \in \{1,2,3 \}) \quad && \in \mathbb{R}^{K \times K}
\\[2mm]
&  \bar{DF}^{-1} := ( \bar{DF}^{-1}_{ab})_{1\le a, b \le 3} \quad && \in \mathbb{R}^{3K\times 3K}\,,
%\\[0mm]
%& \mathbb{U}({\mathbf H}) :=  \mathbb{I}_3 \otimes \text{diag}(e^{-{\hat{U}^0_h}}|_{{\etab}_1} , \cdots,  e^{-{\hat{U}^0_h}}|_{{\etab}_{K}} ) \quad && \in \mathbb{R}^{3K \times 3 K}\,,
\\[2mm]
& \mathbb{W} :=  \mathbb{I}_3 \otimes \text{diag}(w_1, \cdots,  w_{K}) \quad && \in \mathbb{R}^{3K \times 3 K}.
%\\[2mm]
% & \bar{\mathbb{B}}_{\mathbf f}({\mathbf H}) := \text{diag}( \hat{f}({\boldsymbol \eta}_1), \cdots, \hat{f}({\boldsymbol \eta}_K) \quad && \in \mathbb{R}^{3K\times 3K}.
\end{aligned}
\end{equation}

\noindent{\bf Quasi-interpolation}
Here we introduce the local quasi-interpolation projectors detailed in~\cite{3}, which are constructed via the compositions of following one-dimensional local interpolation and histopolation operators:
\begin{equation*}
\begin{aligned}
&\text{interpolation}\quad  I^p f = \sum_{i=0} \lambda_i(f) \hat{N}_i^p, \quad \lambda_i(f) = \sum_{j=0}^{2p-2} \omega^i_j f(x_j^i),\\
&\text{histopolation}\quad H^{p-1} f = \sum_i \tilde{\lambda}_i(f) \hat{D}_i^{p-1}, \quad \tilde{\lambda}_i(f) = \sum_{j=0}^{2p-1}\tilde{\omega}_j^i \int_{\tilde{x}_j^i}^{\tilde{x}^i_{j+1}} f(t) \mathrm{d}t,
\end{aligned}
\end{equation*}
where $p, p-1$ are the degrees of the basis of univariate B-splines $I^p$ and $H^{p-1}$ project onto,  $x^i_j$ and $\tilde{x}^i_{j+1}$ are interpolation and histopolation points related with degrees of freedom $\lambda_i(f)$ and $\tilde{\lambda}_i(f)$, respectively, refer to~\cite{3} for more details and weights $\omega_j$ and $\tilde{\omega}_j$.  
Three dimensional projectors $\Pi_0, {\mathbf \Pi}_1, {\mathbf \Pi}_2, {\Pi}_3$ in commuting diagram~\eqref{diagram} can now be constructed as: 
\begin{equation*}
\begin{aligned}
& \Pi_0 = I^{p_1} \odot I^{p_2} \odot I^{p_3} \\
& {\mathbf \Pi}_1 =  \left(
\begin{matrix}
   H^{p_1-1} \odot I^{p_2} \odot I^{p_3}\\
   I^{p_1} \odot H^{p_2-1} \odot I^{p_3} \\
    I^{p_1} \odot I^{p_2} \odot H^{p_3-1}
      \end{matrix}
\right)  \\
& {\mathbf \Pi}_2 =  \left(
\begin{matrix}
   I^{p_1} \odot H^{p_2-1} \odot H^{p_3-1}\\
   H^{p_1-1} \odot I^{p_2} \odot H^{p_3-1} \\
    H^{p_1-1} \odot H^{p_2-1} \odot I^{p_3}
      \end{matrix}
\right)  \\
& \Pi_3 = H^{p_1-1} \odot H^{p_2-1} \odot H^{p_3-1}, 
\end{aligned}
\end{equation*}
where $\odot$ means "coordinates-wise".
The main feature of the above local projectors is that they are \textit{local}, i.e., the value of $I^pf$ and $H^{p-1}f$ at $\etab$ depends only on the values of $f$ close to $\etab$. This makes the matrices assembled related with local projectors sparse, which helps us construct an efficient solver for sub-step $bv$ in section~\ref{sec:6}. Moreover, corresponding to each projector $\Pi_i$ we introduce the notation $\hat{\Pi}_i: V_i \rightarrow \mathbb{R}^{N_i}$, which gives the finite element coefficients $\lambda_i$ or $\tilde{\lambda}_i,  i = 0, 1, 2, 3$.

\subsection{Bracket discretization}

If we look at the splitting schemes \eqref{stepvv}-\eqref{sys:formIII} there are two obvious ways for phase-space discretization. The first one is to discretize 
the equations directly on the right-hand side with the aim of obtaining skew-symmetric matrices that guarantee energy conservation of each sub-step as~\cite{3}. The second and more systematic way is to discretize the brackets and the Hamiltonians and to derive discrete equations as~\cite{GEMPIC}. The first method works for models without proper brackets (for instance linearized equations), but there can be difficulties when the Hamiltonian functions are not quadratic, as in our case where it contains a logarithmic function. The advantage of the second approach is that anti-symmetric discrete brackets can be obtained directly and thus energy conservations are always guaranteed, no matter how complicated the Hamiltonian functions are.  Also filters could be applied directly in the discretizations of the brackets without destroying  energy conservation properties. 
The discrete Hamiltonian $H $ is given by
\begin{equation}\label{eq:dishhii}
{H} = \frac{1}{2} {\mathbf b}^\top {\mathbb{M}}_1 {\mathbf b} + \frac{1}{2} {\mathbf V}^\top \mathbb{W} {\mathbf V} + \kappa \sum_j q_j\, \hat n^3_h(\etab_j)\ln \left(\hat n^3_h(\etab_j)/\sqrt{g(\etab_j)} \right)\,,
\end{equation}
where the electron thermal energy is approximated by numerical quadrature with weights $q_j$. 
The bracket of our model can now be discretized by substituting $\hat f^3 \to \hat f^3_h$, $\hat\Bb^{1} \to \hat \Bb_h^{1}$ and by inserting the discrete functional derivatives \eqref{eq:weakb} and \eqref{fder:eta}-\eqref{fder:v}. 
For the respective steps of our splitting scheme \eqref{stepvv}-\eqref{sys:formIII} we obtain following discrete sub-brackets:

\subsubsection*{Sub-bracket $\{\hat{\mathcal{P}}, \hat{\mathcal{Q}} \}_{vv}$}
Here we discretize the sub-bracket~\eqref{bras:3} for sub-step~\eqref{stepvv} by substituting directly discrete functional derivatives, and assuming $S$ as a smoothed delta function.
\begin{equation}
\label{eq:rotdis}
\small
\begin{aligned}
& \{ \hat{\mathcal{P}}, \hat{\mathcal{Q}} \}_{vv} \\
& = \int  {\hat f^3} DF^{-\top}\hat\Bb^1 \cdot \nabla_v \frac{\delta \hat{\mathcal{P}}}{\delta \hat f^3} \times \nabla_v \frac{\delta \hat{\mathcal{Q}}}{\delta \hat f^3} \tn{d}{\etab}  \tn{d}{\mathbf v} -\int  \frac{1}{\hat{n}^3}  \hat f^3({\mathbf v}) \hat f^3({\mathbf v}')  DF^{-\top}\hat{\Bb}^1 \cdot \left( \nabla_v \frac{\delta \hat{\mathcal{P}}}{\delta \hat f^3}  \times  \nabla_{v'} \frac{\delta \hat{\mathcal{Q}}}{\delta \hat f^3}  \right) \tn{d}{\etab} \tn{d}{\mathbf v} \tn{d}{\mathbf v}'\,,  \\
& \approx \int \sum_k w_k {S(\etab - \etab_k)}DF^{-\top}\hat\Bb^1_h \cdot \left(\frac{1}{w_k}\frac{\partial P}{\partial {\mathbf v}_k}\right) \times \left(\frac{1}{w_k}\frac{\partial Q}{\partial {\mathbf v}_k}\right) \tn{d}{\etab}\\
& -\int  \frac{1}{\hat{n}^3_h} DF^{-\top} \hat{\Bb}^1_h \cdot \left( \sum_k S(\etab - \etab_k) \frac{\partial P}{\partial {\mathbf v}_k} \times \sum_k S(\etab - \etab_k) \frac{\partial Q}{\partial {\mathbf v}_k} \right)\,  \tn{d}{\etab}\,,\\
& \approx \sum_j q_j \sum_k w_k {S(\etab_j - \etab_k)}DF^{-\top}(\etab_j)\hat\Bb^1_h(\etab_j) \cdot \left(\frac{1}{w_k}\frac{\partial P}{\partial {\mathbf v}_k}\right) \times \left(\frac{1}{w_k}\frac{\partial Q}{\partial {\mathbf v}_k}\right) \\
& -\sum_j q_j  \frac{1}{\hat{n}^3_h(\etab_j)} DF^{-\top}(\etab_j) \hat{\Bb}^1_h(\etab_j) \cdot \left( \sum_k S(\etab_j - \etab_k) \frac{\partial P}{\partial {\mathbf v}_k}  \times \sum_k S(\etab_j - \etab_k) \frac{\partial Q}{\partial {\mathbf v}_k} \right) \\
& =: \nabla_{\mathbf V}P^\top \mathbb{Q}^{vv} \nabla_{\mathbf V}Q\,\\
\end{aligned}
\end{equation}
where $\etab_j$ and $q_j$ are quadrature points and weights, the size of $\mathbb{Q}^{vv}$ (anti-symmetric) is $3K \times 3K$. Numerically, matrix $\mathbb{Q}^{vv}$ is not assembled, and $\nabla_{\mathbf V}P^\top \mathbb{Q}^{vv} \nabla_{\mathbf V}H$ with $P = {\mathbf e}_i \cdot {\mathbf v}_k$ is realized by the directly calculations of the second and third last lines in~\eqref{eq:rotdis}, where   ${\mathbf e}_i$ is the $i$-th unit vector, $i=1,2,3$.
In the calculations in~\eqref{eq:rotdis}, the integrals are approximated by numerical quadratures (scatter-gather method~\cite{scatgath}), and the cost of evaluation of $ \mathbb{Q}^{vv} \mathbb{W} {\mathbf V}$ is $\mathcal{O}(K)$.

\subsubsection*{Sub-bracket $\{\hat{\mathcal{P}}, \hat{\mathcal{Q}} \}_{Bv}$}
Here we discretize the sub-bracket~\eqref{bras:4} for sub-step~\eqref{stepbv}, in which the quasi-interpolation with local property is used to deal with the particle related parts in the bracket, which brings sparse matrices in the numerical schemes detailed in next section.\\
Specifically, by projecting particle parts into finite element space $V_1$ using local projectors, substituting  directly discrete functional derivatives, and assuming $S$ as a smoothed delta function, we have 
\be
\label{eq:disbv3}
 \begin{aligned}
 & \{\hat{\mathcal{P}}, \hat{\mathcal{Q}} \}_{Bv} \\
 & = \int \frac{\hat f^3}{n^3} DF^{-\top}\hat \Bb^1 \cdot \left[ \nabla_v \frac{\delta \hat{\mathcal{P}}}{\delta \hat f^3} \times \left(DF\hat\nabla \times \frac{G}{\sqrt g} \frac{\delta \hat{\mathcal{Q}}}{\delta \hat{\Bb}^1} \right)
  -   \nabla_v \frac{\delta \hat{\mathcal{Q}}}{\delta \hat f^3} \times \left(DF\hat\nabla \times \frac{G}{\sqrt g}\frac{\delta \hat{\mathcal{P}}}{\delta \hat{\Bb}^1}\right)\right] \, \tn{d}{\etab}  \tn{d}{\mathbf v}
  \\
 & = \int \frac{\hat f^3}{n^3} DF^{\top} \nabla_v \frac{\delta \hat{\mathcal{P}}}{\delta \hat{f}^3}
 G^{-1} \left( \hat{\nabla} \times \frac{G}{\sqrt{g}} \frac{\delta \hat{\mathcal{Q}}}{\delta \hat{\mathbf B}^1} \times G^{-1}\hat{\mathbf B}^1 \right) \sqrt{g}\tn{d}{\mathbf v} \tn{d}{\etab}  \\
 & - \int \frac{\hat f^3}{n^3} DF^{\top} \nabla_v \frac{\delta \hat{\mathcal{Q}}}{\delta \hat{f}^3}
 G^{-1} \left( \hat{\nabla} \times \frac{G}{\sqrt{g}} \frac{\delta \hat{\mathcal{P}}}{\delta \hat{\mathbf B}^1} \times G^{-1}\hat{\mathbf B}^1 \right) \sqrt{g}\tn{d}{\mathbf v} \tn{d}{\etab} \\
 & \approx \int \Pi_{1}\left( {DF^{\top}} \sum_k w_k \frac{S(\etab - \etab_k)}{\hat{n}^3_h} \nabla_v \frac{\delta \hat{\mathcal{P}}}{\delta \hat f^3}(\etab, {\mathbf v}_k)\right) \cdot G^{-1} \left[ \left(\mathbb C \mathbb M_1^{-1}\pder{Q}{\bb} \right)^\top \mathbb \Lambda^2 \times G^{-1}\hat \Bb^1_h \right]\sqrt{g} \, \tn{d}{\etab} 
 \\
 &-\int  \Pi_{1} \left( {DF^{\top}}   \sum_k w_k \frac{S(\etab - \etab_k)}{ \hat{n}^3_h} \nabla_v \frac{\delta \hat{\mathcal{Q}}}{\delta \hat f^3}(\etab, {\mathbf v}_k) \right)  \cdot G^{-1} \left[  \left(\mathbb C \mathbb M_1^{-1} \pder{P}{\bb}  \right)^\top \mathbb \Lambda^2 \times G^{-1}\hat \Bb^1_h \right] \sqrt{g}\, \tn{d}{\etab}  \tn{d}{\mathbf v}\\
 & \approx \nabla_{\mathbf V}P^\top \mathbb{L} \mathbb{Q}^{bv} \mathbb{C} \mathbb{M}_1^{-1} \nabla_{\mathbf b}Q - \nabla_{\mathbf V}Q^\top \mathbb{L} \mathbb{Q}^{bv} \mathbb{C} \mathbb{M}_1^{-1} \nabla_{\mathbf b}P,\\
 & \approx \nabla_{\mathbf V}P^\top \mathbb{L} \mathbb{T}^\top \mathbb{Q}^{bv} \mathbb{C} \mathbb{M}_1^{-1} \nabla_{\mathbf b}Q - \nabla_{\mathbf V}Q^\top \mathbb{L} \mathbb{T}^\top \mathbb{Q}^{bv} \mathbb{C} \mathbb{M}_1^{-1} \nabla_{\mathbf b}P.
\end{aligned}
\ee
where we choose $\Pi_1$ as the quasi-interpolation local projector, 
and $\mathbb{L} \in \mathbb{R}^{3K \times N_1}$ with components as 
\begin{equation}\label{eq:LLex}
\begin{aligned}
\mathbb{L}_{k, i} =  \left(\Pi_1\left(\frac{S(\etab - \etab_k)}{\hat{n}^3_h(\etab)} {DF^\top(\etab)}{\mathbf e}_1  \right)\right)_i,\\
\mathbb{L}_{k+K, i} = \left(\Pi_1\left(\frac{S(\etab - \etab_k)}{\hat{n}^3_h(\etab)}{DF^\top(\etab)}{\mathbf e}_2 \right)\right)_i,\\
\mathbb{L}_{k+2K, i} =  \left(\Pi_1\left(\frac{S(\etab - \etab_k)}{\hat{n}^3_h(\etab)}{DF^\top(\etab)}{\mathbf e}_3  \right)\right)_i.
\end{aligned}
\end{equation}
where $1 \le k \le K$, $1 \le i \le N_1$, ${\mathbf e}_s$ is the $s$-th unit vector in $\mathbb{R}^3$, $s = 1 ,2, 3$. We also include filters for noise reduction in our discrete brackets. Bilinear filters or variants thereof~\cite{filter} will be used. ${\mathbb{T}}$ denotes the matrix of filters, and we introduced the matrix $\mathbb Q^{bv} \in \mathbb R^{N_1 \times N_2}$ with entries
$$
 \mathbb Q^{bv}_{ij} := \sum_l q_l \Lambda^1_i(\etab_l) \cdot G^{-1}(\etab_l) \left( \Lambda^2_j(\etab_l) \times G^{-1}(\etab_l)\hat{\mathbf B}_h^1(\etab_l)   \right) \sqrt{g(\etab_l)},
$$
where $q_l$ and $\etab_l$ are quadrature weights and points. 

%\begin{remark}
%To use local projectors, the smoothed delta function is necessary to use. In the above calculations, we made the following approximation,
%$$
%w_k S(\etab - \etab_k) \nabla_v \frac{\delta \hat \cF}{\delta \hat f^3}(\etab, {\mathbf v}_k) \approx S(\etab - \etab_k) \frac{\partial F}{\partial {\mathbf v}_k}.
%$$
%%which gives the same discretizations of ${\mathbf b}$ by directly writing the equations~\eqref{stepvv} in a weak formulation with local projectors and smoothed delta functions. 
%\end{remark}

%From the "strong" bracket \eqref{bra:2:strong} we obtain
%\be
%\begin{aligned}
% \tn{strong:}\quad  \{ \hat\cF, \hat\cH \}_{Bv} &\approx \int \frac{\hat f^3_h}{\hat D^3_h} \left[ \nabla_v \frac{\delta \hat \cF}{\delta \hat f^3}  \cdot DF^{-\top} \left(\hat\nabla \times \pder{H^\top}{\bb} \mathbb M_2^{-1} \frac{G}{\sqrt g} \mathbb \Lambda^2 \right) \times \hat \Bb^2_h \right. 
% \\[1mm]
% &  \left. \qquad\qquad -  \nabla_v \frac{\delta \hat \cG}{\delta \hat f^3} \cdot DF^{-\top}  \left(\hat\nabla \times \pder{F_2^\top}{\bb} \mathbb M_2^{-1} \frac{G}{\sqrt g} \mathbb \Lambda^2 \right) \times \hat \Bb^2_h \right] \tn{d}{\etab}  \tn{d}{\mathbf v}\,,  \nonumber
% \end{aligned}
% \ee

\subsubsection*{Sub-bracket $\{\hat{\mathcal{P}}, \hat{\mathcal{Q}} \}_{BB}$}
For the sub-bracket~\eqref{bras:5} corresponding to~\eqref{stepbb}, we have the following discretization.
\begin{equation}
\label{eq:disbb}
 \begin{aligned}
& \{ \hat{\mathcal{P}}, \hat{\mathcal{Q}} \}_{BB} = -\int \frac{1}{\hat{n}^3}  DF^{-\top}\hat{\Bb}^1 \cdot \left(DF \hat \nabla \times \frac{G}{\sqrt g} \frac{\delta \hat{\mathcal{P}}}{\delta \hat{\Bb}^1}\right) \times \left(DF \hat\nabla \times \frac{G}{\sqrt g} \frac{\delta \hat{\mathcal{Q}}}{\delta \hat{\Bb}^1}\right) \, \tn{d}{\etab}\,,  
     \\[1mm]
    & \approx -\int \frac{1}{\hat{n}^3}  DF^{-\top}\hat{\Bb}^1_h \cdot \left(DF (\mathbb{\Lambda}^2)^\top \mathbb{C} \mathbb{M}_1^{-1} \frac{\partial P}{\partial {\mathbf b}} \right) \times \left(DF (\mathbb{\Lambda}^2)^\top \mathbb{C} \mathbb{M}_1^{-1}  \frac{\partial Q}{\partial {\mathbf b}} \right) \, \tn{d}{\etab},\\[1mm]
    & = -\int \frac{1}{\hat{n}^3}  G^{-\top}\hat{\Bb}^1_h \cdot \left( (\mathbb{\Lambda}^2)^\top \mathbb{C} \mathbb{M}_1^{-1} \frac{\partial P}{\partial {\mathbf b}} \right) \times \left(  (\mathbb{\Lambda}^2)^\top \mathbb{C} \mathbb{M}_1^{-1} \frac{\partial Q}{\partial {\mathbf b}} \right) \sqrt{g}\, \tn{d}{\etab},\\
    & \approx  \left(\frac{\partial P}{\partial {\mathbf b}} \right)^\top \mathbb{M}_1^{-1}   \mathbb{C}^\top \mathbb{Q}^{bb}({\mathbf b})  \mathbb{C} \mathbb{M}_1^{-1}   \frac{\partial Q}{\partial {\mathbf b}},
 \end{aligned}
\end{equation}
where $ \mathbb{Q}^{bb}_{ij}({\mathbf b}) = - \sum_l q_l  \frac{\sqrt{g}(\etab_l)}{\hat{n}^3_h(\etab_l)} G^{-\top}(\etab_l)\hat{\Bb}^1_h(\etab_l) \cdot \left( \Lambda^2_i(\etab_l) \times  \Lambda^2_j(\etab_l)  \right)$, $\etab_l$ and $q_l$ are quadrature points and weights.

\begin{remark}
The discretization derived from \eqref{eq:disbb} is 
$
\dot{\mathbf b} = \mathbb{M}_1^{-1}   \mathbb{C}^\top \mathbb{Q}^{bb}  \mathbb{C} {\mathbf b},
$ 
which is the same as the scheme obtained by discretizing the equations~\eqref{stepbb} in a weak formulation and assuming $\hat{B}^1_h \in V_1$.
\end{remark}

\subsubsection*{Sub-bracket $\{\hat{\mathcal{P}}, \hat{\mathcal{Q}} \}_{xv}$}
For the sub-bracket $\{\hat{\mathcal{P}}, \hat{\mathcal{Q}} \}_{xv}$~\eqref{bras:1}, we get the following discretization assuming that $S$ is a delta function.
\begin{equation}
\label{eq:disxviii}
 \begin{aligned}
 & \{\hat{\mathcal{P}}, \hat{\mathcal{Q}} \}_{xv} \\
  &= \int \hat f^3  \left( DF^{-\top} \hat\nabla_\eta \frac{\delta \hat{\mathcal{P}}}{\delta \hat f^3}  \cdot \nabla_v \frac{\delta \hat{\mathcal{Q}}}{\delta \hat f^3} 
 - \nabla_v\frac{\delta \hat{\mathcal{P}}}{\delta \hat f^3} \cdot DF^{-\top}\hat \nabla_\eta  \frac{\delta \hat{\mathcal{Q}}}{\delta \hat f^3} \right) \tn{d}{\etab}  \tn{d}{\mathbf v} \,,\\[1mm]
 & = \int \sum_{k=1}^{K} w_k \delta({\etab} - {\etab}_k ) \delta({\mathbf v} - {\mathbf v}_k )   \left( DF^{-\top} \hat\nabla_\eta \frac{\delta \hat{\mathcal{P}}}{\delta \hat f^3}  \cdot \nabla_v \frac{\delta \hat{\mathcal{Q}}}{\delta \hat f^3}  -  \nabla_v\frac{\delta \hat{\mathcal{P}}}{\delta \hat f^3} \cdot DF^{-\top}\hat \nabla_\eta  \frac{\delta \hat{\mathcal{Q}}}{\delta \hat f^3}  \right)  \tn{d}{\etab}  \tn{d}{\mathbf v}, \\[1mm]
 & =  \sum_{k=1}^{K} \left( \frac{1}{w_k} DF^{-\top}(\etab_k) \frac{\partial P}{\partial {\etab}_k }  \cdot \frac{\partial Q}{\partial {\mathbf v}_k} - \frac{1}{w_k} DF^{-\top}(\etab_k) \frac{\partial Q}{\partial {\etab}_k }  \cdot \frac{\partial P}{\partial {\mathbf v}_k}  \right). 
 \end{aligned}
\end{equation}

\section{Time discretizations}\label{sec:6}
Time discretizations are done using the Poisson splitting introduced in section~\ref{sec:2}, i.e., via splitting the bracket into several parts. In this section, we present the details of all the sub-steps. The $\delta f$ method can be derived similarly as in~\cite{3} by decomposing the distribution function as the sum of equilibrium and non-equilibrium parts.
For sub-steps $vv, bv, bb$, mid-point rules are used in time, and for sub-step $xv$, Hamiltonian splitting methods and discrete gradient methods are used. For an unknown $a$, $a^n$ and $a^{n+1}$ mean the numerical approximations of $a$ at $n$-th and $n+1$-th time step, respectively. Moreover, $a^{n+\frac{1}{2}} = \frac{a^n + a^{n+1}}{2}$.

\subsection{Sub-step $vv$}\label{sec:vv}
From sub-bracket~\eqref{eq:rotdis} and by setting $P  = {\mathbf v}_k, 1 \le k \le K$ and $Q=H$, we get the following equations for particles,
\begin{equation}\label{eq:vvlocal}
\dot{\mathbf V} = \mathbb{Q}^{vv} \mathbb{W}{\mathbf V}.
\end{equation}
From~\eqref{stepvv}, we know that every particle is related with all the particles through the current term, then the cost of usual methods would be $\mathcal{O}(K^2)$. The scatter-gather method~\cite{scatgath} is used in~\eqref{eq:rotdis}, which reduces the cost of evaluation of right hand of the above equation~\eqref{eq:vvlocal} to $\mathcal{O}(K)$. Mid-point rule is used for time discretization, 
$$
\frac{{\mathbf V}^{n+1} - {\mathbf V}^{n} }{\Delta t} = \mathbb{Q}^{vv} \mathbb{W} \frac{{\mathbf V}^{n} + {\mathbf V}^{n} }{2}.
$$
which conserves the energy. Note that momentum $\sum_k w_k {\mathbf v}_k$ is also conserved by this sub-step. In fact, we have, for $i = 1, 2, 3$ 
\begin{equation*}
\begin{aligned}
& \left(\sum_k w_k {\mathbf v}_k^{n+1}\right)_i - \left(\sum_k w_k {\mathbf v}_k^{n}\right)_i \\
& = \Delta t \sum_j q_j \sum_k DF^{-\top}(\etab_j)\hat\Bb^1_h(\etab_j) \cdot {\mathbf e}_i  \times   {S(\etab_j - \etab_k)} w_k {\mathbf v}_k^{n+\frac{1}{2}} \\
& - \Delta t \sum_j q_j  DF^{-\top}(\etab_j) \hat{\Bb}^1_h(\etab_j) \cdot \left( \underbrace{\frac{\sum_k w_k S(\etab_j - \etab_k) }{ \hat{n}^3_h(\etab_j)}}_{=1}{\mathbf e}_i \times \sum_k S(\etab_j - \etab_k) w_k {\mathbf v}_k^{n+\frac{1}{2}}   \right) = 0,
\end{aligned}
\end{equation*}
where ${\mathbf e}_i$ is the $i$-th unit vector.
%As the particle number is relatively large,  implicit mid-point rule (iterations are needed) is costly to use. In order to get energy conservation efficiently, we use a second order Runge--Kutta method combined with a linear scaling method proposed in \cite{Jiang1}.
%The final scheme is explicit, thus efficient to use. 

\subsection{Sub-step $bv$}
From sub-bracket~\eqref{eq:disbv3}, and by setting $P = {\mathbf v}_k, {\mathbf b}, 1 \le k \le K$ and $Q=H$, we get 
\begin{equation}\label{eq:bvdiss}
\dot{\mathbf V} = \mathbb{L} \mathbb{T}^\top \mathbb{Q}^{bv} \mathbb{C} {\mathbf b}, \quad \dot{\mathbf b} = - \mathbb{M}_1^{-1} \mathbb{C}^\top  (\mathbb{Q}^{bv})^\top \mathbb{T} \mathbb{L}^\top   \mathbb{W} {\mathbf V},
\end{equation}
for which midpoint rule is used to conserve energy, i.e., 
\begin{equation*}
\begin{aligned}
\frac{{\mathbf V}^{n+1} -  {\mathbf V}^{n}  }{\Delta t} &=  \mathbb{L} \mathbb{T}^\top \mathbb{Q}^{bv} \mathbb{C} \frac{{\mathbf b}^{n+1} +  {\mathbf b}^{n}  }{2},\\
 \frac{{\mathbf b}^{n+1} -  {\mathbf b}^{n}  }{\Delta t} &= -\mathbb{M}_1^{-1}  \mathbb{C}^\top (\mathbb{Q}^{bv})^\top \mathbb{T} \mathbb{L}^\top   \mathbb{W} \frac{{\mathbf V}^{n+1} +  {\mathbf V}^{n}  }{2}.
\end{aligned}
\end{equation*}
After substituting the equation about ${\mathbf V}$ into the equation about ${\mathbf b}$, we get the following scheme
\begin{equation*}
\begin{aligned}
& \left(\mathbb{M}_1 + \frac{\Delta t^2 }{4} \mathbb{C}^\top (\mathbb{Q}^{bv})^\top \mathbb{T} \mathbb{L}^\top  \mathbb{W}  \mathbb{L}  \mathbb{T}^\top  \mathbb{Q}^{bv} \mathbb{C}  \right) {\mathbf b}^{n+1} \\
& = \left(\mathbb{M}_1 - \frac{\Delta t^2 }{4}\mathbb{C}^\top (\mathbb{Q}^{bv})^\top \mathbb{T} \mathbb{L}^\top  \mathbb{W}  \mathbb{L} \mathbb{T}^\top \mathbb{Q}^{bv} \mathbb{C}  \right) {\mathbf b}^{n} - \Delta t \mathbb{C}^\top (\mathbb{Q}^{bv})^\top \mathbb{T} \mathbb{L}^\top \mathbb{W}  {\mathbf V}^n,
\end{aligned}
\end{equation*}
which is solved using a fixed point iteration method and conjugate gradient method  with a pre-conditioner of $\mathbb{M}_1$ (by fast Fourier transformation).
\begin{remark}\label{bvremark}
The divergence free property of the magnetic field, i.e.,  $\nabla \cdot {\mathbf B} = 0$ holds weakly. In fact, we have 
$$
\mathbb{G}^\top \mathbb{M}_1 {
\mathbf b}^{n+1} = \mathbb{G}^\top \mathbb{M}_1 {
\mathbf b}^{n} - \Delta t   \mathbb{G}^\top \mathbb{C}^\top (\mathbb{Q}^{bv})^\top \mathbb{T} \mathbb{L}^\top   \mathbb{W} \frac{{\mathbf V}^{n+1} +  {\mathbf V}^{n}  }{2} = 0, 
$$
which uses the property that $ \mathbb{G}^\top \mathbb{C}^\top = 0$.
The momentum $\sum_k w_k {\mathbf v}_k$ is not conserved by this sub-step. In fact, the conservation of momentum is obtained with the condition that $\nabla \cdot {\mathbf B} = 0$ in a strong sense. 
\end{remark}

\begin{remark}
The sparse property of $\mathbb{L}$ comes from using local quasi-interpolation projectors, which makes the matrix $\mathbb{L}^\top \mathbb{W} \mathbb{L}$ also sparse. As the positions of particles are unchanged during this sub-step, the matrix  $\mathbb{L}^\top \mathbb{W} \mathbb{L}$ only needs to be assembled before the start of the fixed point iterations. 
\end{remark}

\subsection{Sub-step $bb$}
From sub-bracket~\eqref{eq:disbb}, by setting  $P = {\mathbf b}$ and $Q=H$, we get 
\begin{equation}\label{eq:bbdis}
\dot{\mathbf b} =   \mathbb{M}_1^{-1} \mathbb{C}^\top \mathbb{Q}^{bb}({\mathbf b})   \mathbb{C} {\mathbf b},
\end{equation}
As the Hamiltonian is quadratic about ${\mathbf b}$, the mid-point rule is used for \eqref{eq:bbdis} to conserve energy, i.e., 
$$
\frac{{\mathbf b}^{n+1} -  {\mathbf b}^{n} }{\Delta t} =  \mathbb{M}_1^{-1} \mathbb{C}^\top \mathbb{Q}^{bb}({\mathbf b}^{n+\frac{1}{2}})   \mathbb{C}  {\mathbf b}^{n+\frac{1}{2}}.
$$
By multiplying $\mathbb{M}_1$ on both sides, we have the following equivalent formulation 
\begin{equation*}
\begin{aligned}
 \left( \mathbb{M}_1 - \frac{\Delta t}{2} \mathbb{C}^\top \mathbb{Q}^{bb}({\mathbf b}^{n+\frac{1}{2}}) \mathbb{C} \right) {\mathbf b}^{n+1}  = \mathbb{M}_1 {\mathbf b}^n + \frac{\Delta t}{2}  \mathbb{C}^\top \mathbb{Q}^{bb}({\mathbf b}^{n+\frac{1}{2}}) \mathbb{C}   {\mathbf b}^{n}.
\end{aligned}
\end{equation*}
The fixed point iteration scheme is, (from $i$-th iteration to  $(i+1)$-th iteration),
\begin{equation*}
\begin{aligned}
 \left( \mathbb{M}_1 - \frac{\Delta t}{2} \mathbb{C}^\top \mathbb{Q}^{bb}({\mathbf b}^{n+\frac{1}{2}}_i) \mathbb{C} \right) {\mathbf b}^{n+1}_{i+1}  = \mathbb{M}_1 {\mathbf b}^n + \frac{\Delta t}{2}  \mathbb{C}^\top \mathbb{Q}^{bb}({\mathbf b}^{n+\frac{1}{2}}_i) \mathbb{C}   {\mathbf b}^{n}, \ {\mathbf b}^{n+\frac{1}{2}}_i = \frac{ {\mathbf b}^n + {\mathbf b}^n_{i} }{2}.
\end{aligned}
\end{equation*}
In each iteration, the above linear system is solved using GMRES method with an incomplete LU decomposition or inverse fast Fourier transformation of the mass matrix (of periodic boundary condition and cartesian coordinates) as a pre-conditioner.
As ${\mathbf b}$ is changing every fixed point iteration, assembling the matrix (every iteration) $\mathbb{Q}^{bb}$ is costly, which instead can be avoided by doing matrix-free matrix vector products. 
%defining a linear operator: 
%\begin{equation*}
%\begin{aligned}
% Op: \mathbb{R}^{N_2} \rightarrow  \mathbb{R}^{N_2},  \, & {\xb} \rightarrow {\mathbf y} =   \mathbb{Q}^{bb}({D}, {\mathbf b}^{n+\frac{1}{2}}_k) {\xb}, \, y_i =   \int  {\mathbf \Lambda}^2_i({\etab})   \cdot  \left(  G^{-1}  ({\mathbf \Lambda}^1)^\top {\mathbf b}_k^{n+\frac{1}{2}}  \times ({\mathbf \Lambda}^2)^\top {\xb} \right) \frac{1}{{D}}   \tn{d}{\etab }.
%\end{aligned}
%\end{equation*}
Similar to remark~\ref{bvremark}, we can prove that $\nabla \cdot {\mathbf B} = 0 $ weakly, i.e., $\mathbb{G}^\top \mathbb{M}_1 {
\mathbf b}^{n+1} = \mathbb{G}^\top \mathbb{M}_1 {
\mathbf b}^{n}$.

\subsection{Sub-step $xv$}
\label{subsec:xvn}

From sub-brackets~\eqref{eq:disxviii}, by setting $P  = \etab_k, {\mathbf v}_k, 1 \le k \le K$ and $Q=H$, we get 
\begin{equation}
\label{eq:nlnqn}
\begin{aligned}
& \left(
\begin{matrix}
    \dot{\mathbf{H}} \\
     \dot{\mathbf{V}}
      \end{matrix}
\right)  \!=\!\! 
\left(
\begin{matrix}
    {\mathbf{0}} \!\! &   \bar{DF}^{-1}({\mathbf H})  {\mathbb{W}}^{-1} \\
     -   {\mathbb{W}}^{-1}\bar{DF}^{-\top}({\mathbf H}) \!\! & {\mathbf 0}
           \end{matrix}
\right) 
\left(
\begin{matrix}
    \nabla_{\mathbf{H}}  H  \\
     \nabla_{\mathbf{V}}  H = \mathbb{W} {\mathbf V}
      \end{matrix}
\right).
\end{aligned}
\end{equation}
Corresponding Hamiltonian is 
\begin{equation*}
\begin{aligned}
H & =  \frac{1}{2} {\mathbf V}^\top \mathbb{W} {\mathbf V} + \frac{1}{2}{\mathbf b}^\top \mathbb{M}_1 {\mathbf b} + \kappa  \sum_j \sum_{k=1}^K q_j {w_{k}} S({\etab}_j - {\etab}_{k}) \ln \left(\sum_{k=1}^K \frac{w_{k}}{\sqrt{g({\etab}_j)}}  S({\etab}_j - {\etab}_{k}) \right).
\end{aligned} 
\end{equation*}
%\begin{remark}
%The above scheme can be proved as a good approximation for the original equation with the help of integration by parts,
%\begin{equation}
%\begin{aligned}
%\nabla_{{\etab}_k} H & = \kappa  \int \left( \ln  \left(\sum_{k'=1}^K \frac{w_{k'}}{\sqrt{g}}  S({\etab} - {\etab}_{k'}) \right) + 1 \right) (-\frac{w_k}{\sqrt{g}} ) \nabla S({\etab} - {\etab}_{k}) \sqrt{g}  \tn{d}{\etab},\\
%& = \kappa  \int  \ln  \left(\sum_{k'=1}^K \frac{w_{k'}}{\sqrt{g}}  S({\etab} - {\etab}_{k'}) \right)  (-\frac{w_k}{\sqrt{g}} )  \nabla S({\etab} - {\etab}_{k}) \sqrt{g}  \tn{d}{\etab},\\ 
%&  = \kappa  w_k  \int \frac{\nabla \left( \sum_{k'=1}^K \frac{w_{k'}}{\sqrt{g}}  S({\etab} - {\etab}_{k'})\right) }{\sum_{k'=1}^K \frac{w_{k'}}{\sqrt{g}}  S({\etab} - {\etab}_{k'}) } S({\etab} - {\etab}_{k})    \tn{d}{\etab} \\
%& \approx \kappa  w_k \frac{\nabla \left( \sum_{k'=1}^K \frac{w_{k'}}{\sqrt{g}}  S({\etab} - {\etab}_{k'}) \right) = \nabla n({\etab})}{\sum_{k'=1}^K \frac{w_{k'}}{\sqrt{g}}  S({\etab} - {\etab}_{k'}) =n({\etab})  }|_{{\etab}_k}\\
%\end{aligned} 
%\end{equation}
%\end{remark}
\noindent{\bf Cartesian coordinates case}
System~\eqref{eq:nlnqn} is a  Hamiltonian system, for which we use Hamiltonian splitting method and get the following two sub-steps
\begin{equation}
\label{eq:IIhamiltoniansplit}
\begin{aligned}
\dot{\mathbf V} &=  -   {\mathbb{W}}^{-1}\bar{DF}^{-\top}  \nabla_{\mathbf{H}}  H, \\
\dot{\mathbf H} &=  \bar{DF}^{-1}  {\mathbf{V}},
\end{aligned}
\end{equation}
which both can be solved explicitly.
Although Hamiltonian splitting makes the energy not conserved exactly by the whole scheme, we will see in the numerical experiments section that energy and momentum are still conserved with high accuracy over long time. Also as it is explicit, it is an important choice for the large scale simulations in cartesian coordinates. \\
\noindent{\bf Curvilinear coordinates case}
For this case, sub-step $\dot{\mathbf H} =  \bar{DF}^{-1}({\mathbf H})  {\mathbf{V}}$ given by Hamiltonian splitting methods can not be solved explicitly. The following midpoint discrete gradient method~\cite{Gonzalez} can be used to conserve the energy exactly.
\begin{equation*}
\begin{aligned}
& \frac{{\mathbf H}^{n+1} - {\mathbf H}^{n}}{\Delta t} =   {\mathbb{W}}^{-1} \bar{DF}^{-1}({\mathbf H}^{n+\frac{1}{2}}) \bar{\nabla}_{\mathbf{V}}  H,\\
& \frac{{\mathbf V}^{n+1} - {\mathbf V}^{n}}{\Delta t} = -   {\mathbb{W}}^{-1}\bar{DF}^{-\top}({\mathbf H}^{n+\frac{1}{2}})  \bar{\nabla}_{\mathbf{H}}  H,
\end{aligned}
\end{equation*}
where 
\begin{equation*}
\begin{aligned}
& \bar{\nabla}_{\mathbf H} H = \nabla_{\mathbf H} H \left(\frac{{\mathbf H}^n + {\mathbf H}^{n+1}  }{2}\right) + f_{disg}\left( {\mathbf H}^{n+1} - {\mathbf H}^n \right),\\
& \bar{\nabla}_{\mathbf V} H = \nabla_{\mathbf V} H \left(\frac{{\mathbf V}^n + {\mathbf V}^{n+1}  }{2}\right) + f_{disg}\left( {\mathbf V}^{n+1} - {\mathbf V}^n \right).
\end{aligned}
\end{equation*} 
where 
$$
f_{disg} = \frac{H({\mathbf H}^{n+1}, {\mathbf V}^{n+1}) - H({\mathbf H}^n, {\mathbf V}^{n}) - \nabla H (\frac{{\mathbf H}^n + {\mathbf H}^{n+1}  }{2}, \frac{{\mathbf V}^n + {\mathbf V}^{n+1}  }{2}) \cdot (  ({\mathbf H}^{n+1} - {\mathbf H}^n)^\top,   ({\mathbf V}^{n+1} - {\mathbf V}^n)^\top)^\top    }{ |{\mathbf H}^{n+1} - {\mathbf H}^n|^2 +  |{\mathbf V}^{n+1} - {\mathbf V}^n|^2 }.
$$

\begin{remark}
In the above discrete gradient methods, the degrees of shape functions should be at least 2 in each direction, which guarantees that $\nabla_{\mathbf H}H$ are continuous functions about ${\mathbf H}$, and is vital to guarantee the convergences of the iterations in discrete gradient methods. 
\end{remark}

In summary, to get the complete numerical schemes, we run the above four integrators successively, i.e., 
\begin{align}
& \Phi_{\Delta t} := \Phi^{xv, ha}_{\Delta t}   \circ \Phi^{vv}_{\Delta t}   \circ \Phi^{bb}_{\Delta t}   \circ \Phi^{bv}_{\Delta t}, \label{eq:LieI},\\
& \Phi_{\Delta t} := \Phi^{xv, dis}_{\Delta t}   \circ \Phi^{vv}_{\Delta t}   \circ \Phi^{bb}_{\Delta t}   \circ \Phi^{bv}_{\Delta t}, \label{eq:LieII}
\end{align}
which are first order Lie-Trotter splitting methods~\cite{Trotter}, where $\Phi^{vv}_{\Delta t}, \Phi^{bv}_{\Delta t}, \Phi^{bb}_{\Delta t}$ are solution maps for sub-steps $vv, bv, bb$, respectively, $\Phi^{xv, ha}_{\Delta t}$ and $\Phi^{xv, dis}_{\Delta t} $ are solution maps for sub-step $xv$ using Hamiltonian splitting methods and discrete gradient methods, respectively. Higher order schemes can be obtained using composition methods~\cite{HLW}.

\section{Comparisons}\label{sec:7}

\subsection{Conserved Properties}
Here we summarize all the properties conserved by the above numerical methods~\eqref{eq:LieI}-\eqref{eq:LieII}. Positivity of electron density is conserved by both schemes, which is represented as the sum of finite number of particles. As for energy conservation, it could be obtained by~\eqref{eq:LieII} with discrete gradient methods used for sub-step-\eqref{eq:nlnqn}, and numerically we will see that energy is also conserved with a very high precision by~\eqref{eq:LieI} with explicit Hamiltonian splitting methods used in sub-step-\eqref{eq:nlnqn} in cartesian coordinate case. Both schemes~\eqref{eq:LieI} and~\eqref{eq:LieII} conserve $ \nabla \cdot {\mathbf B} = 0$ (weakly). The quasi-neutrality relation is built in the model equations.  
While momentum $\sum_k w_k {\mathbf v}_k$ is only conserved by sub-steps \eqref{eq:bbdis} and \eqref{eq:vvlocal}, numerically we always see that momentum is conserved with high accuracy.

\subsection{Performance}
The implementations of the above methods of two formulations are done in the framework of python package STRUPHY~\cite{3} with OpenMP/MPI hybrid parallelisation. Particle related parts are parallelized using MPI and OpenMP, and linear systems in sub-steps are solved by iterative methods on one node currently (with up to 72 threads). 

Here we mainly consider the performances of the particle related parts as \cite{3}.
 The performances of two above constructed methods are compared by running with the same parameters in cartesian coordinates. The time used by particle related sub-steps $xv$, $vv$, $bv$ scales well with thread number. For sub-step $xv$, discrete gradient methods $\Phi^{xv, dis}_{\Delta t}$ in sub-step $xv$ are usually more costly than explicit Hamiltonian splittings  $\Phi^{xv, ha}_{\Delta t} $. Moreover, solving the linear systems in sub-steps $bb$ and $bv$ with multi-nodes with MPI is already on the way to improve the efficiency further.

%\begin{table}[hbt]
%\begin{tabular}{ccccccc}
%\hline
%\hline
%Openmp threads  & substep $xvn$ & substep $vv$ & substep $bb$ &  substep $bv$ \\
%10  & 1.15 &  11.14 & 0.47  &  5.45    \\
%20  & 0.68 & 5.57   & 0.24  &  2.85  \\
%40  & 0.40 & 2.88 & 0.13 & 1.61     \\
%\hline
%\end{tabular}
%\caption{Formulation I ($fn{\Bb}$), one time step, $Np = 60000$, $[200, 2, 2]$ cells.}
%\vskip 0.2cm 
%\begin{tabular}{cccccccc}
%\hline
%\hline
%MPI   & substep 1 & substep 2 & substep 3 & substep 4 (II) & substep 4 (III) & substep 5 &  substep 6  \\
%1  & 1.658 & 2.824 & 5.754 & 17.881  & & 0.137 & 0.143 \\
%2  & 1.603&  2.375&5.572  & 17.852   & &  0.094 &  0.079 \\
%4  & 2.052 & 2.481  &  8.214 & 22.289 & & 0.119  & 0.113 \\
%\hline
%\end{tabular}
%\caption{ $fU{\Bb}$, 100 time steps, $Np = 10000$, $[20, 2, 2]$ cells.}
%\vskip 0.2cm 
\begin{table}[hbt]
\begin{tabular}{ccccccc}
\hline
\hline
threads  & $\frac{\text{substep}\ xv}{\text{iteration number}}~\eqref{eq:LieII}$ & substep $xv$~\eqref{eq:LieI} & $\frac{\text{substep}\ vv}{\text{iteration number}} $ &  substep $bv$  \\
  9 & 4.93 & 2.46& 2.94   & 9.60   \\
  18 & 2.43 & 1.23& 1.46   & 4.95   \\
 36 &1.29 & 0.61 & 1.01  & 3.04 \\
 72 & 0.62 & 0.31   & 0.38 & 1.70  \\
\hline
\end{tabular}
\caption{One time step, particle number is $60000$, cells number is $[200, 3, 3]$.}
\end{table}

\section{Numerical experiments}\label{sec:8}
In this section, some numerical tests are done to validate the implementations: finite grid instability, Landau damping, R wave, Bernstein waves, mirror instability, and parametric instability. 
The following two kinds of mapping (cartesian, colella) $\mathbf{F}$ from logical domain to physical domain  are used in the numerical simulations
\begin{equation}
\label{eq:Colella}
\begin{aligned}
& \text{Cartesian}: \quad \mathbf{F}: \hat{\Omega} \rightarrow \Omega, \quad \etab \rightarrow \left(
\begin{matrix}
  L_x \eta_1 \\
  L_y  \eta_2\\
 L_z \eta_3
       \end{matrix}
\right) = {\mathbf x},\\
& \text{Colella}: \quad \mathbf{F}: \hat{\Omega} \rightarrow \Omega, \quad \etab \rightarrow \left(
\begin{matrix}
  L_x (\eta_1 + \alpha \sin(2\pi \eta_1)\sin(2\pi \eta_2))\\
  L_y (\eta_2 + \alpha \sin(2\pi \eta_2)\sin(2\pi \eta_3))\\
 L_z \eta_3
       \end{matrix}
\right) = {\mathbf x}.
\end{aligned}
\end{equation}
We check the conservations properties mentioned above. Moreover, dispersion relations of some numerical tests are checked numerically. Iteration tolerance is set to be $10^{-12}$. As the density appears in the denominator in some terms of this model, when the density is very low, there would be some numerical instabilities. We follow the strategy in~\cite{current2D}  to set the relevant terms as zero when the density is 0 or below a minimum threshold, which is set to be $10^{-3}$ in the following simulations. And the support parameters $h_1, h_2, h_3$ of shape functions in~\eqref{eq:shape_support} are chosen as the cell sizes in each space direction if we do not mention specifically. In the following, we refer to the schemes~\eqref{eq:LieI}-\eqref{eq:LieII} by Hamiltonian splitting method and discrete gradient method, respectively. Periodic boundary conditions are considered.

%We test with two kinds of mapping, i.e., Cartesian mesh and Colella mesh,
%\begin{equation}
%\begin{aligned}
%&\text{Cartesian}: {\mathbf F} : \hat{\Omega} \rightarrow \Omega, \quad \etab \rightarrow 
%\left(
%\begin{matrix}
% L_x \eta_1 \\
%L_y \eta_2 \\
%L_z \eta_3
%           \end{matrix}
%\right)  = {\mathbf x},\\
%&\text{Colella}: {\mathbf F} : \hat{\Omega} \rightarrow \Omega, \quad \etab \rightarrow 
%\left(
%\begin{matrix}
% L_x (\eta_1 + \alpha\sin(2\pi\eta_1)\sin(2\pi \eta_2) ) \\
% L_x (\eta_12+ \alpha\sin(2\pi\eta_2)\sin(2\pi \eta_3) ) \\
%L_z \eta_3
%           \end{matrix}
%\right)  = {\mathbf x}.
%\end{aligned}
%\end{equation}
%Here physical domain is $[L_x, L_y, L_z]$, and $0 \le \alpha < \frac{1}{2\pi}$.

\subsection{Finite grid instability}
\begin{figure}[htbp]
\center{\includegraphics[scale=0.38]{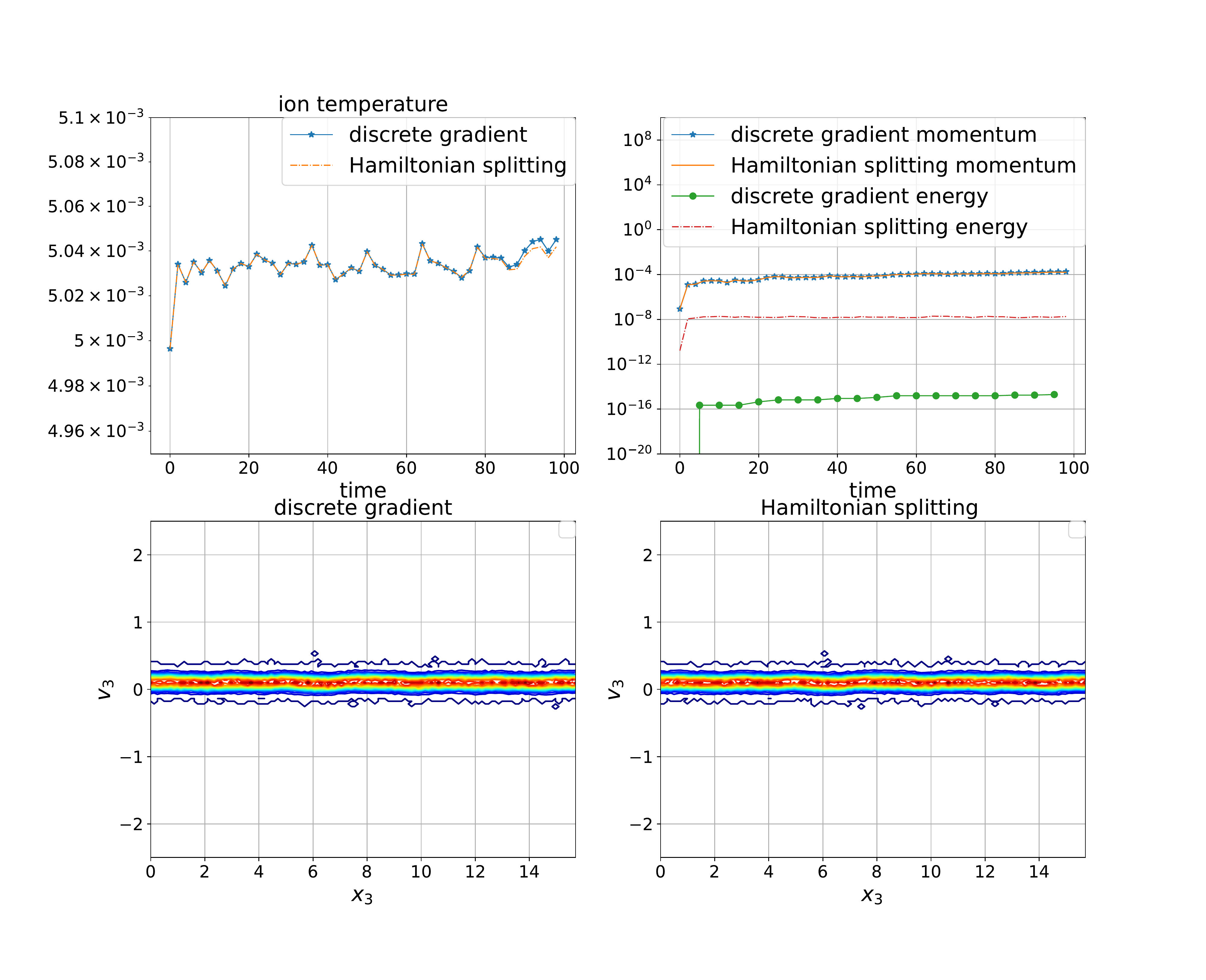}}
\caption{{\bf Finite grid instability.} Time evolution of the temperature of ions, momentum errors (the third component), and the $x_3-v_3$ contour plots of schemes~\eqref{eq:LieI}-\eqref{eq:LieII} at $t=100$. }
\label{fig:finite}
\end{figure}

As~\cite{chacon1}, we check the numerical behaviors of the discretizations~\eqref{eq:LieI}-\eqref{eq:LieII} on the finite grid instability. The initial conditions (equilibrium of the system) are 
\begin{equation*}
\begin{aligned}
 {\mathbf B} = {\mathbf 0}, \ f =  \frac{1}{\pi^{\frac{3}{2}} v_T^{\frac{3}{2}}}  e^{-\frac{|v_x|^2}{v_T^2} - \frac{|v_y|^2}{v_T^2}  - \frac{|v_z - 0.1|^2}{v_T^2} }, \,  \kappa = 1, \, v_T = 0.1.
\end{aligned}
\end{equation*}
Cartesian mapping is used for this numerical experiment.
Computational parameters are: 
grid number $[4, 4, 32]$, domain parameters $[L_x,  L_y,  L_z] = [1, 1, 5\pi]$, degrees of B-splines $[1, 1, 3]$, time step size $\Delta t = 0.005$, degrees of shape functions $[2, 2, 2]$, quadrature points in each cell $[2, 2, 4]$, and total particle number $2\times10^5$.  As pointed out in~\cite{Rambo}, usual PIC methods have the finite grid instability, i..e the growth of temperature of ions,  due to the errors induced by the deposition and interpolation procedures. We run the simulations with the above initial conditions and the numerical methods~\eqref{eq:LieI}-\eqref{eq:LieII}, and the results are presented in Fig.~\ref{fig:finite}. We can see that both schemes capture the equilibrium numerically. As ${\mathbf B} = 0$ during the numerical simulations, we are in fact only solving sub-step $xv$.  As the Hamiltonian splitting method (Strang splitting) for sub-step $xv$ in~\eqref{eq:LieI} is a symplectic method, it has superior long numerical behaviors, although energy is not conserved exactly (with a relative error about $10^{-8}$), and the temperature of ions is almost unchanged even after 1 million time steps (not shown). The discrete gradient method~\eqref{eq:LieII} conserves the energy exactly,  also the time evolution of temperature of ions is almost the same as that given by Hamiltonian splitting method during $t \in [0, 100]$. Note that the momentum is not conserved by~\eqref{eq:LieI}-\eqref{eq:LieII} exactly, but its error is at the level of $10^{-4}$, and has no obvious and quick growth with time. 
Also from the other figures of Fig.~\ref{fig:finite}, we can see the $x_3-v_3$ contour plots of the distribution function at $t = 100$ given by these two schemes, the particles are still close to a very thin Maxwellian, due to the good conservation of energy. 

\subsection{Ion Landau damping}
\begin{figure}[htbp]
\center{
\includegraphics[scale=0.4]{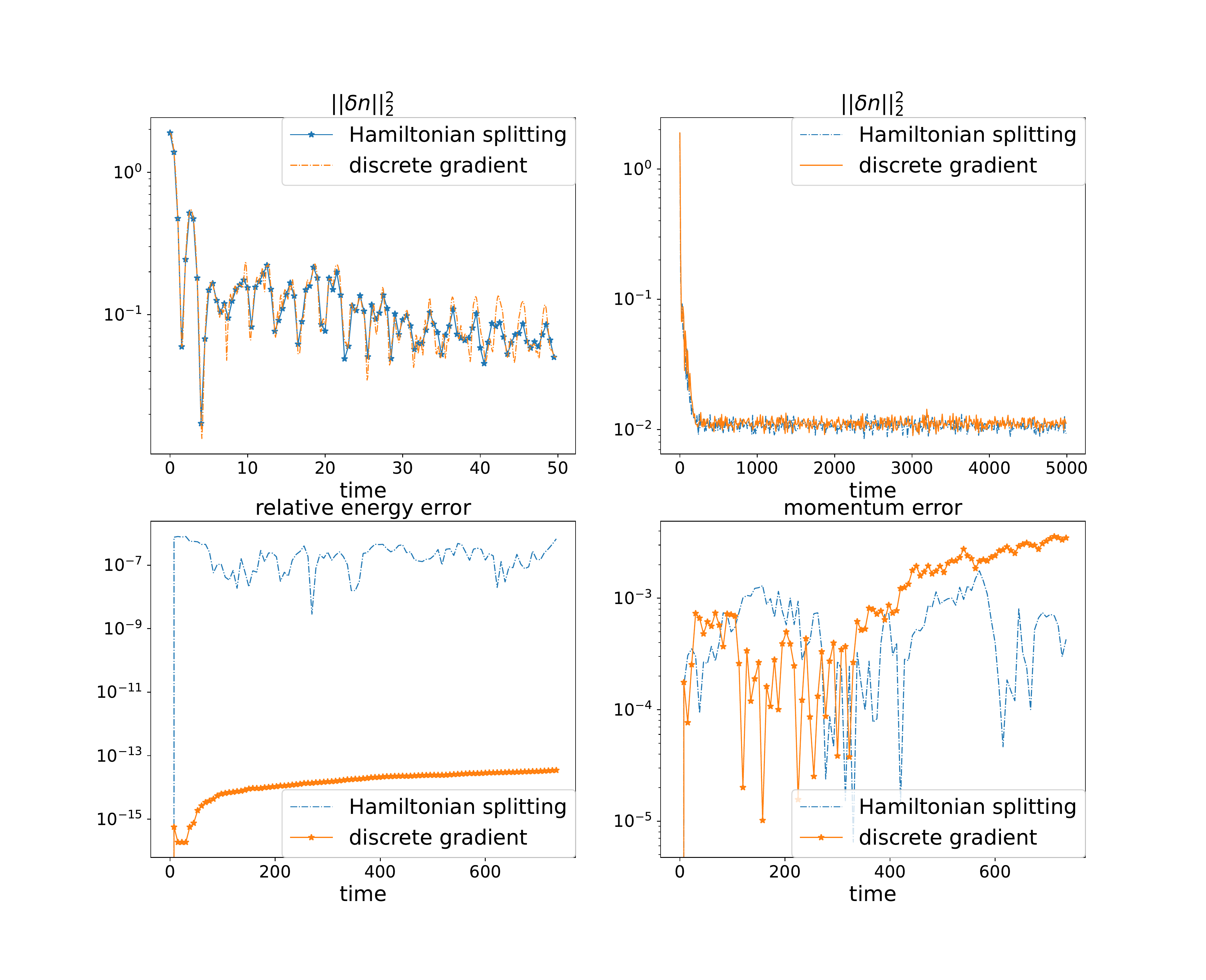}
%\subfigure[]{\includegraphics[scale=0.49]{fnquasi-neutral.eps}}
%\subfigure[]{\includegraphics[scale=0.49]{fULandau_damping.eps}}
%\subfigure[]{\includegraphics[scale=0.49]{fULandau_energyerror.eps}}
%\subfigure[]{\includegraphics[scale=0.49]{fnbLandau_energyerror.eps}}
}
\caption{{\bf Landau damping} Time evolution of $\int |\delta n|^2 \sqrt{g}\mathrm{d}{\etab}$, $\delta n = n - 1$, relative energy errors, and momentum errors (the third component).  }
\label{fig:Landau_all}
\end{figure}

Then we simulate ion Landau damping by a quasi-1D dimensional simulation without background magnetic field. Initial conditions are:
\begin{equation*}
\begin{aligned}
 {\Bb} = {\mathbf 0},  \ f =  \frac{1}{\pi^{\frac{3}{2}} v_T^{\frac{3}{2}}} (1 + 0.5\cos(0.4z)) e^{-\frac{|v_x|^2}{v_T^2} - \frac{|v_y|^2}{v_T^2}  - \frac{|v_z|^2}{v_T^2} }, \,  \kappa = 6.25. 
\end{aligned}
\end{equation*}
Cartesian mapping is used for this numerical experiment.
Computational parameters are:
grid number $[3, 3, 50]$, domain parameters $[L_x,  L_y,  L_z] = [1, 1, 5\pi]$, time step size $\Delta t = 0.005$, final computation time $5\times 10^3$, $v_T = 1.4142$, total particle number $5 \times 10^4$, degrees of shape functions $[2, 2, 2]$, quadrature points in each cell $[2, 2, 2]$, and the support parameters in~\eqref{eq:shape_support} $h_i = \Delta \eta_i, i = 1, 2, 3$. 
See the numerical results in Fig.~\ref{fig:Landau_all}. Both schemes~\eqref{eq:LieI}-\eqref{eq:LieII} give similar results, and $\int |\delta n|^2 \sqrt{g}\mathrm{d}{\etab}$ decays exponentially initially, then grows exponentially with time, and it saturates at a small value ($10^{-2}$) after  $t=200$.  The relative energy error of scheme~\eqref{eq:LieI} is about $10^{-6}$ and has no obvious growth with time, and the scheme~\eqref{eq:LieII} conserves energy at the level of $10^{-15}$. The momentum is conserved by both schemes at high precisions with an error at the level of $10^{-3}$. 
 
\subsection{Parallel electromagnetic wave: R mode} 
\begin{figure}[htbp!]
\center{
\subfigure[]{\includegraphics[scale=0.49]{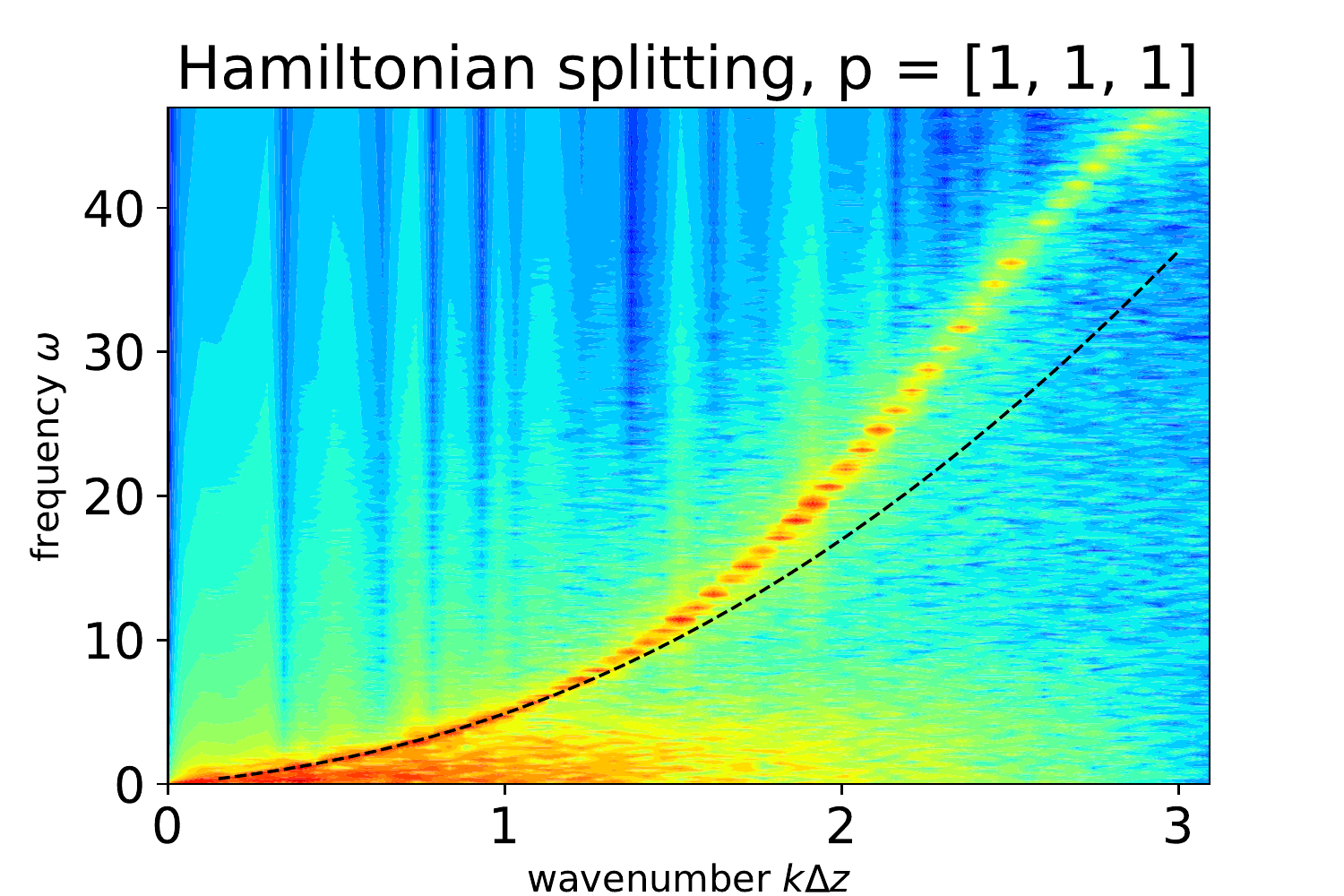}}
\subfigure[]{\includegraphics[scale=0.49]{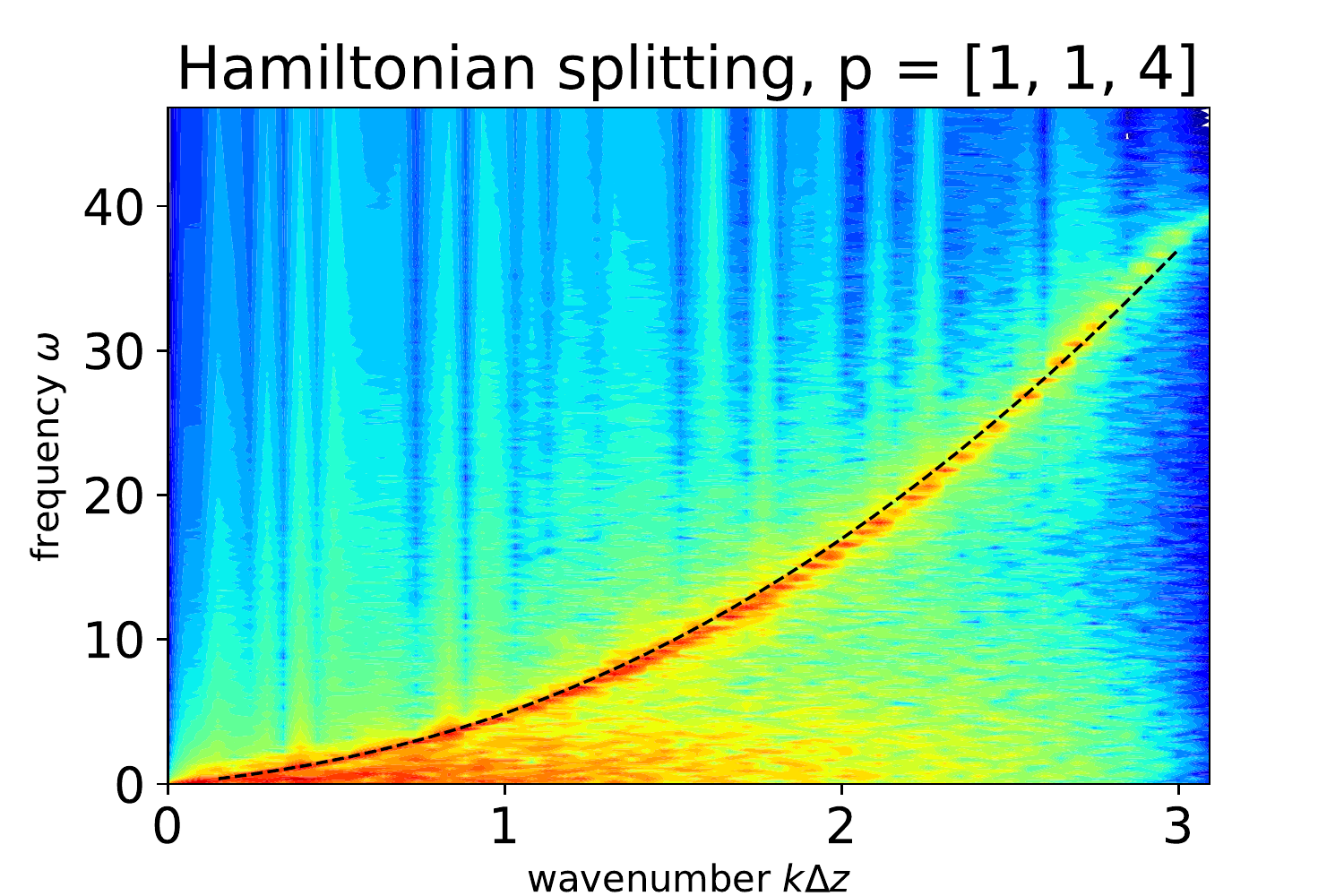}}
%\subfigure[]{\includegraphics[scale=0.49]{fU1_50000.eps}}
%\subfigure[]{\includegraphics[scale=0.49]{fU4_50000.eps}}
\subfigure[]{\includegraphics[scale=0.49]{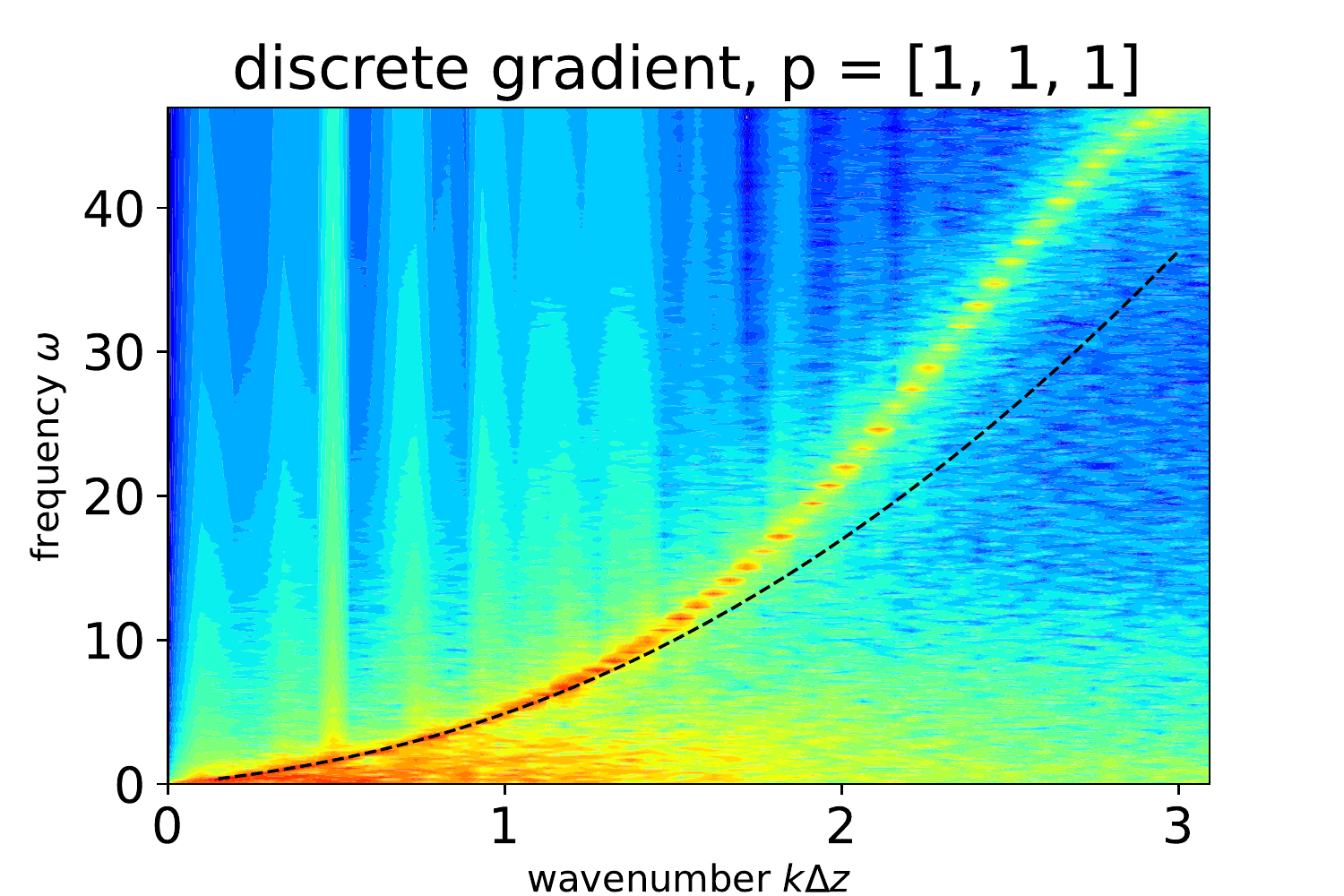}}
\subfigure[]{\includegraphics[scale=0.49]{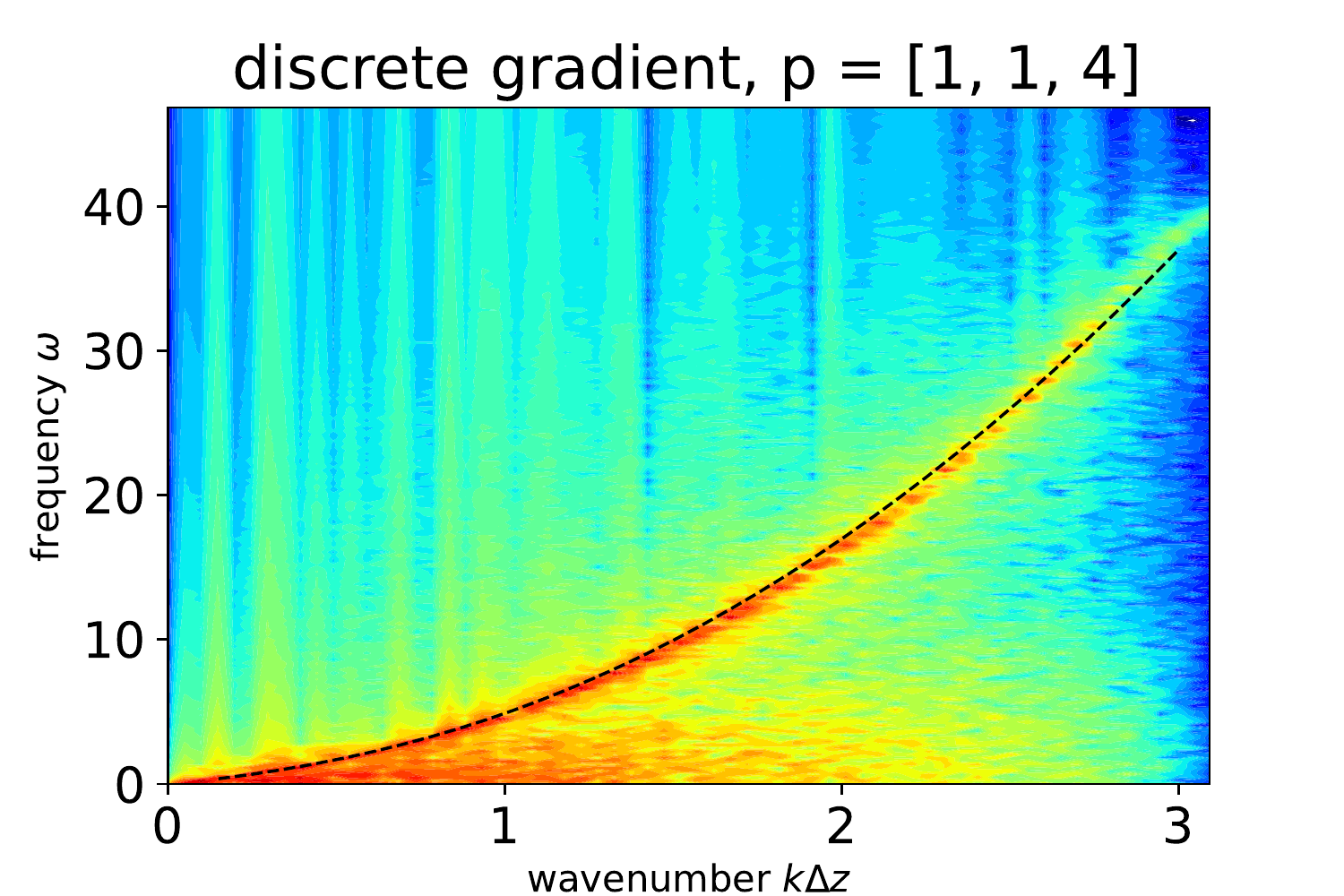}}
}
\caption{{\bf R mode} (a) numerical dispersion of scheme~\eqref{eq:LieI} with $p=[1, 1, 1]$;  (b) numerical dispersion of scheme~\eqref{eq:LieI} with $p=[1, 1, 4]$; (c) numerical dispersion of scheme~\eqref{eq:LieII} with $p=[1, 1, 1]$; (d) numerical dispersion of scheme~\eqref{eq:LieII} with $p=[1, 1, 4]$. }
\label{fig:Rmode1_all}
\end{figure}

\begin{figure}[htbp]
\center{
\includegraphics[scale=0.4]{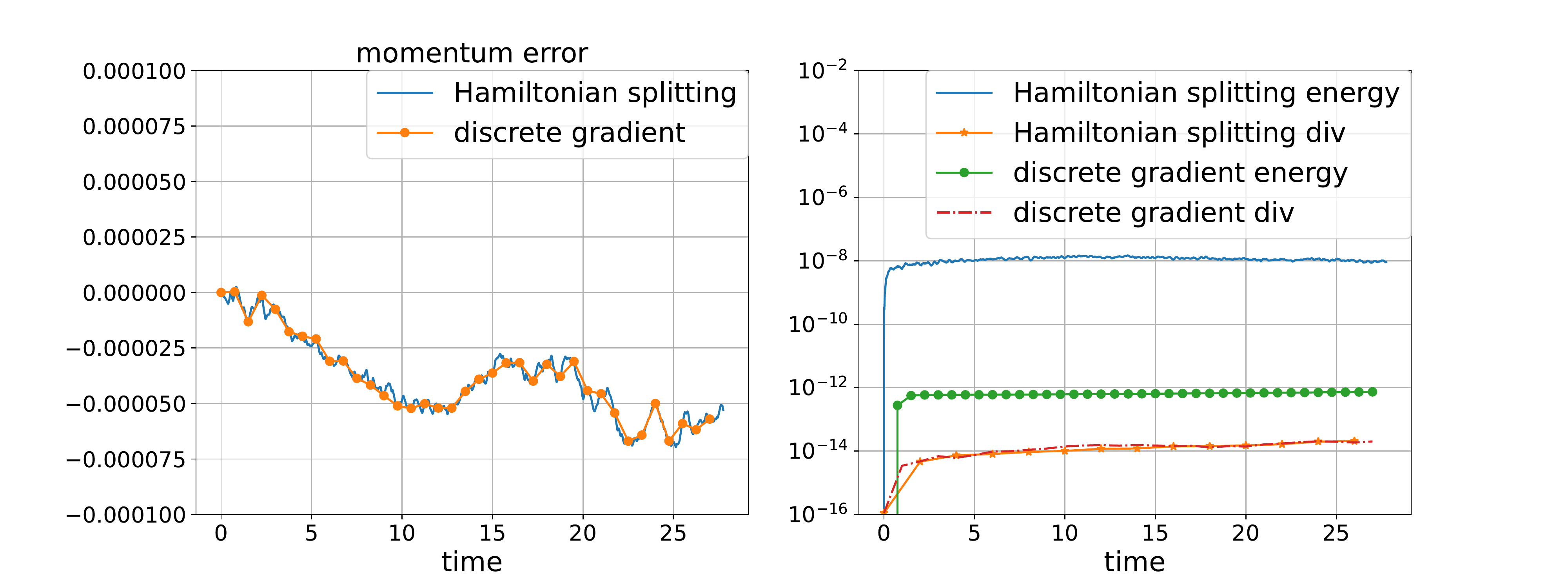}
}
\caption{{\bf R mode}: time evolutions of momentum errors (the third component), relative energy errors, and weak divergence of ${\mathbf B}$ of schemes~\eqref{eq:LieI}-\eqref{eq:LieII} with $p=[1, 1, 4]$.}
\label{fig:Rmode2_all}
\end{figure}

In this part, we check a parallel propagating wave--R wave by a quasi-1D (mostly along one direction, averaging over 2 cells in the transverse directions) simulation. Background magnetic field is along $z$ direction. We do not apply any perturbation for the system other than the noise of the PIC method due to the reduced number of macro-particles. Specifically, initial conditions we used are:
\begin{equation*}
\begin{aligned}
{\Bb} = (0, 0, 1), \, \kappa  = 1, \, f =  \frac{1}{\pi^{\frac{3}{2}} }e^{-{|v_x|^2} - {|v_y|^2}  - {|v_z|^2} }.
\end{aligned}
\end{equation*}
A cartesian mapping is used for this numerical experiment.
Computational parameters are: 
grid number $[3, 3, 128]$, domain parameters $[L_x,  L_y,  L_z] = [1, 1, 64]$, time step size $\Delta t = 0.005$, final computation time $40$, total particle number $5\times 10^{5}$, degrees of shape functions $[2, 2, 2]$, and quadrature points in each cell $[2, 2, 4]$. Full-f method is used for this simulation. The numerical results of the dispersion relation of the R wave are presented in Fig.~\ref{fig:Rmode1_all}. The black dashed lines correspond to the analytical dispersion relations given by python package HYDRO proposed in \cite{disp}, when $k \ll 1$, $\omega \propto k$, when $k \gg 1$, $\omega \propto k^2$. We can see that our numerical results are in good agreement with the analytical results when the degrees of the basis functions are $[1, 1, 4]$, even when the wave number $k$ is larger than Nyquist frequency $\frac{\pi}{2\Delta z}$. However, when the degrees of B-spline are lower, for example, $p=[1,1,1]$, we can see that there are more obvious errors in the numerical dispersion relations when $k$ is larger than Nyquist frequency.
The results of conservation properties of schemes~\eqref{eq:LieI}-\eqref{eq:LieII} are presented in  Fig.~\ref{fig:Rmode2_all}. We can see that the relative energy error of schemes~\eqref{eq:LieI}-\eqref{eq:LieII} are at the level of $10^{-8}$, and $10^{-12}$, respectively. The errors of the third component of momentum are also displayed, and we can see that there is no linear or quick growth, and the momentum errors of both schemes~\eqref{eq:LieI}-\eqref{eq:LieII} are about $10^{-5}$. As we proved that $\nabla \cdot {\mathbf B} = 0$ in the weak sense, the numerical errors are at the level of $10^{-14}$. 

%We also plot the electron density at some specific time, $t = ...$ using the three numerical methods constructed above. We can see that ...

%\begin{figure}[htbp]
%\center{
%\subfigure[]{\includegraphics[scale=0.45]{R_momentum_high.eps}}
%}
%\caption{{\bf R mode}: time evolutions of the third component of momentum error with $p=[1, 1, 4]$.}
%\label{fig:R_momentum_all}
%\end{figure}

\subsection{Perpendicular wave: ion Bernstein wave}
\begin{figure}[htbp!]
\center{
{\includegraphics[scale=0.44]{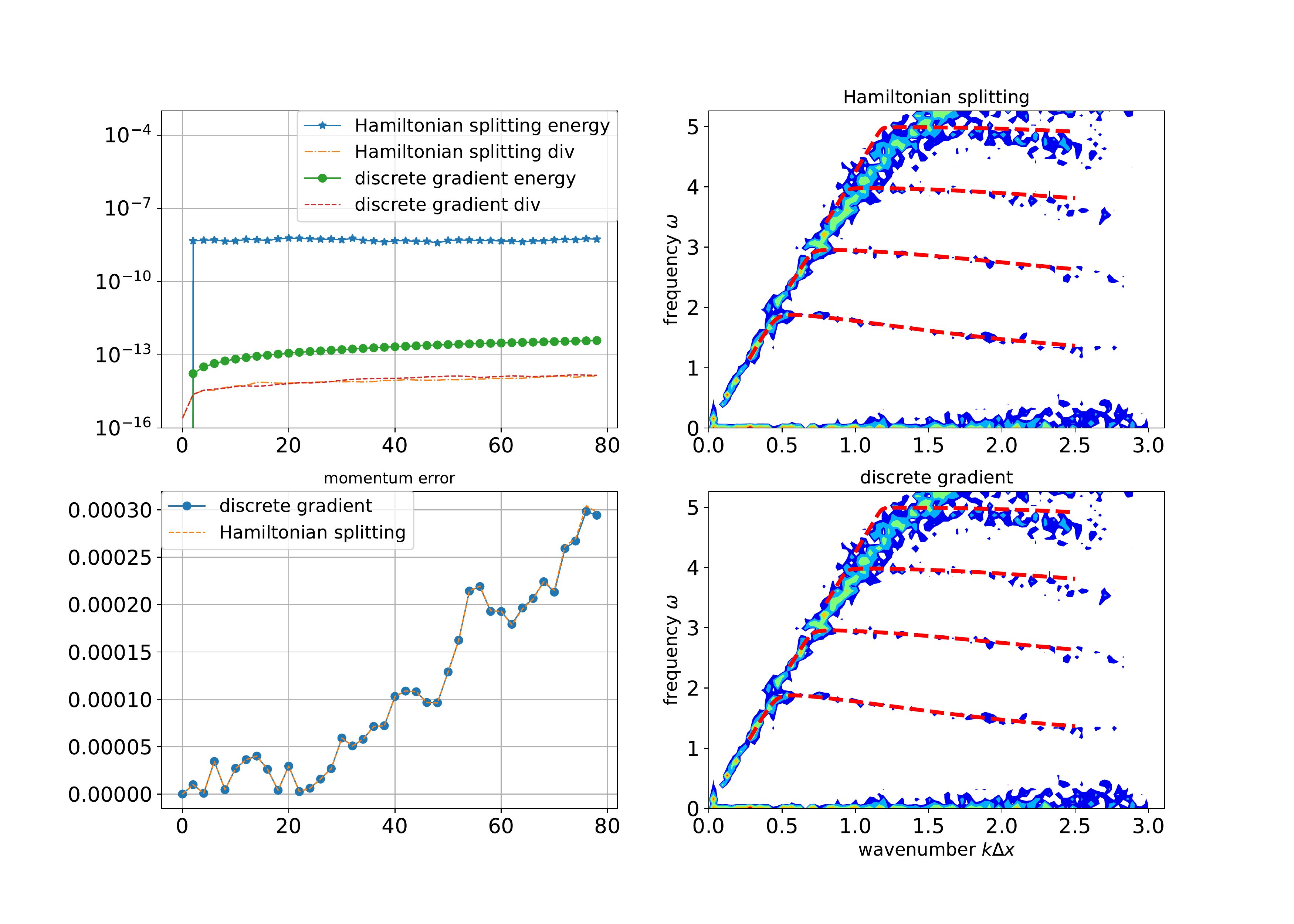}}
}
\caption{{\bf Bernstein wave 1D.} Time evolutions of relative energy errors, weak divergence of magnetic field, momentum errors, and dispersion relations of Bernstein waves schemes~\eqref{eq:LieI}-\eqref{eq:LieII}.}
\label{fig:ber_all}
\end{figure}

Next we check a wave nearly perpendicular to a magnetic field--Bernstein wave by a one dimensional simulation, which only exists in a plasma with finite temperature. 
In order to excite these waves, we initialize a quasi-1D thermal plasma along the $x$ direction. No initial perturbation is added except the noise of the PIC method. Specifically, initial conditions are:
\begin{equation*}
\begin{aligned}
 {\Bb} = (0, 0, 1), \ f =  \frac{1}{\pi^{\frac{3}{2}} v_T^{\frac{3}{2}}} e^{-\frac{|v_x|^2}{v_T^2} - \frac{|v_y|^2}{v_T^2}  - \frac{|v_z|^2}{v_T^2} }, \, \kappa  = 0.09. 
\end{aligned}
\end{equation*}
Cartesian mapping is used for this numerical experiment.
Computational parameters are: 
grid number $[200, 3, 3]$, domain parameters $[L_x,  L_y,  L_z] = [50, 1, 1]$, time step size $\Delta t = 0.01$, $v_T = 0.2121$, final computation time $80$. Particle number $6 \times 10^4$, degrees of polynomials $[4,1,1]$, degrees of shape functions $[2, 2, 2]$. A full-f method is used in these two simulations. See the time evolution of relative energy errors, momentum errors, and dispersion relations of Bernstein waves in Fig.~\ref{fig:ber_all}, in which red dashed lines are analytical dispersion relations of Bernstein waves obtained via HYDRO code~\cite{disp}. We can see that both methods~\eqref{eq:LieI}-\eqref{eq:LieII} can give rather good dispersion relations numerically. 
And both schemes give almost the same momentum errors. The relative energy errors of ~\eqref{eq:LieI} is  at the level of $10^{-8}$, and the relative energy error of~\eqref{eq:LieII} is at the level of $10^{-13}$. 
Note that there are some errors for the dispersion relation when the frequency is about 5 in Fig.~\ref{fig:ber_all}, which is due to the smoothing effects of shape functions, especially when the degree of shape function is high. But as pointed out in~\cite{martin}, this phenomenon can be improved by using more grid points (higher resolution) in space. And the numerical errors of $\nabla \cdot {\mathbf B} = 0$ (weak error) are at the level of $10^{-14}$. 
% \alertred{When more efficient local projection (smoothed delta function is a must to use) is used in sub-step 4, incorrect numerical dispersion relation is obtained because of the smoothing effects of smoothed delta functions~ \cite{martin}. To obtain accurate numerical dispersion relation, better resolution (more grid points) is required, which is demonstrated in Fig.~\ref{fig:ber_smooth} with $f{\mathbf B}$ formulation.  }

%\begin{figure}[htbp!]
%\center{
%{\includegraphics[scale=0.49]{B_momentum_high.eps}}
%}
%\caption{{\bf Bernstein wave 1D.} Time evolutions of the error of the third momentum for Formulation I and II.}
%\label{fig:ber_all2}
%\end{figure}

%\begin{figure}[htbp!]
%\center{
%\subfigure[]{\includegraphics[scale=0.49]{local_Bern.eps}}
%\subfigure[]{\includegraphics[scale=0.49]{local_Bern_old.eps}}
%}
%\caption{{\bf Results of Local projection of $f{\mathbf B}$ formulation}: (a) 3D domain is $[0,25]\times  [0,1] \times [0, 1]$, cell number is $[200, 2, 2]$, support of shape function is $[0.5, 1, 1]$; (b) 3D domain is $[0,50]\times  [0,1] \times [0, 1]$, cell number is $[200, 2, 2]$, support of shape function is $[1, 1, 1]$. In both cases, shape functions degree is $[3, 1, 1]$ (in $x, y, z$ direction, respectively).}
%\label{fig:ber_smooth}
%\end{figure}

\subsection{Mirror instability}
\begin{figure}[htbp]
\center{
{\includegraphics[scale=0.42]{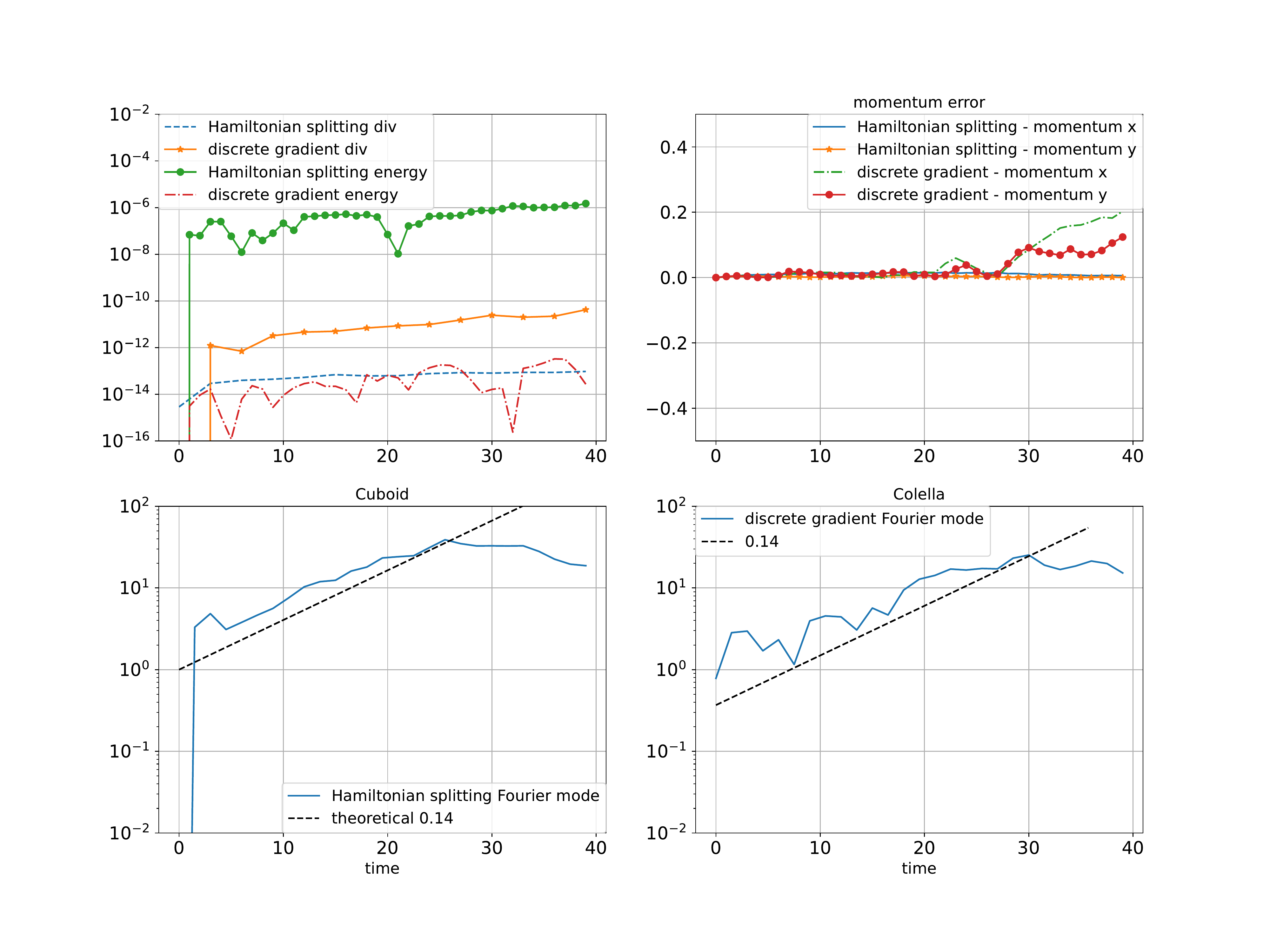}}
}
\caption{{\bf Mirror instability 2D}. Time evolution of the relative energy errors, divergence of magnetic field, momentum errors (the first and second component), and Fourier modes $\sqrt{|\hat{B}_x({k_x, k_y})|^2 + \hat{B}_y({k_x, k_y})|^2}$, $k_x = \frac{2\pi}{L_x}, k_y = \frac{2\pi}{L_y}$.}
\label{fig:fire2Dmirrior}
\end{figure}

Previous numerical tests are all quasi-1D, next we validate our code by a quasi-2D test, the mirror instability. 
Initial conditions are 
\begin{equation*}
\begin{aligned}
 &{\Bb} = (1, 0, 0), \,  f =  \frac{1}{\pi^{\frac{3}{2}} (v_{T_x} v_{T_y} v_{T_z}     )^{\frac{1}{2}}} e^{-\frac{|v_x|^2}{v_{Tx}^2} - \frac{|v_y|^2}{v_{Ty}^2}  - \frac{|v_z|^2}{v_{Tz}^2} }, \, \kappa = 1.
\end{aligned}
\end{equation*}
Cartesian and colella (with $\alpha = 0.05$ in~\eqref{eq:Colella}) mappings are used for this numerical experiment. The scheme~\eqref{eq:LieI} with explicit Hamiltonian splitting~\eqref{eq:IIhamiltoniansplit} (Strang splitting) is used for the cartesian mapping case, and the scheme~\eqref{eq:LieII} with discrete gradient method is used for colella mapping case.
Computational parameters are: 
grid number $[20, 20, 2]$, domain parameters $[L_x,  L_y,  L_z] = [20, 20, 2]$, $k_x =\frac{2\pi}{L_x}$, $k_y = \frac{2\pi}{L_y}$, time step size $\Delta t = 0.01$, $v_{Ty} = v_{Tz} = \sqrt{15}$, $v_{Tx}=\sqrt{6}$, and total particle number $2.5 \times 10^5$. 
We run the simulations for  the two cases: cartesian mapping case and colella mapping case. 
See the numerical results of the two schemes in Fig.~\ref{fig:fire2Dmirrior}, we can see that accurate numerical instability rates ($\gamma = 0.14 $ obtained by HYDRO~\cite{disp}) of Fourier mode $\sqrt{|\hat{B}_x({k_x, k_y})|^2 + \hat{B}_y({k_x, k_y})|^2}$, $k_x = \frac{2\pi}{L_x}, k_y = \frac{2\pi}{L_y}$ are given for the cartesian coordinate case, also growth rates do not change too much in the colella mapping case. Relative energy error is at the level of $10^{-7}$ for the cartesian mapping case, and energy is conserved in the iteration tolerance level in the colella mapping case. There are some growths of the momentum errors of scheme~\eqref{eq:LieII} after $t = 25$. And the numerical errors of $\nabla \cdot {\mathbf B} = 0$ (weak error) are at the level of $10^{-14}$.

%\begin{figure}[htbp]
%\center{
%\subfigure[]{\includegraphics[scale=0.49]{fUB_2D_Bx.png}}
%\subfigure[]{\includegraphics[scale=0.49]{fUB_2D_By.png}}
%\subfigure[]{\includegraphics[scale=0.49]{fUB_2D_Bz.png}}
%\subfigure[]{\includegraphics[scale=0.49]{fUB_2D_density.png}}
%}
%\caption{{\bf Mirror instability 2D, $\delta f$, $fU{\Bb}$ formulation, $5 \times 10^6$ particles, $\alpha = 10^{-3}$ }: (a) Time evolution of the Fourier mode $ { |B_x(k_x, k_y)|}$; (b)  Time evolution of the Fourier mode $ { |B_y(k_x, k_y)|}$; (c)  Time evolution of the Fourier mode $ { |B_z(k_x, k_y)|}$;  (d) Time evolution of the Fourier mode $ { |e^U(k_x, k_y)|}$; $k_x = \frac{2\pi}{20}$, $k_y = \frac{2\pi}{20}$.}
%\label{fig:mirror2DfUB}
%\end{figure}

%\begin{figure}[htbp]
%\center{
%\subfigure[]{\includegraphics[scale=0.49]{ypert_Bx_2D.png}}
%\subfigure[]{\includegraphics[scale=0.49]{ypert_By_2D.png}}
%\subfigure[]{\includegraphics[scale=0.49]{ypert_Bz_2D.png}}
%\subfigure[]{\includegraphics[scale=0.49]{ypert_density_2D.png}}}
%\caption{{\bf Mirror instability 2D, $\delta f$, $fU{\Bb}$ formulation, $2 \times 10^6$ particles, only add perturbation of $B_y = 10^{-3}\cos(\frac{\pi}{10}x)  \cos(\frac{\pi}{10}y)$ }: (a) Time evolution of the Fourier mode $ { |B_x(k_x, k_y)|}$; (b)  Time evolution of the Fourier mode $ { |B_y(k_x, k_y)|}$; (c)  Time evolution of the Fourier mode $ { |B_z(k_x, k_y)|}$;  (d) Time evolution of the Fourier mode $ { |e^U(k_x, k_y)|}$; $k_x = \frac{2\pi}{20}$, $k_y = \frac{2\pi}{20}$.}
%\label{fig:mirror2DfUBy}
%\end{figure}

%\begin{figure}[htbp]
%\center{
%\subfigure[]{\includegraphics[scale=0.49]{fB_Bx_2D.png}}
%\subfigure[]{\includegraphics[scale=0.49]{fB_By_2D.png}}
%\subfigure[]{\includegraphics[scale=0.49]{fB_Bz_2D.png}}
%\subfigure[]{\includegraphics[scale=0.49]{fB_density_2D.png}}
%}
%\caption{{\bf Mirror instability 2D, $\delta f$, $f{\Bb}$ formulation, $1.5\times 10^5$ particles, $\alpha = 10^{-3}$ }: (a) Time evolution of the Fourier mode $ { |B_x(k_x, k_y)|}$; (b)  Time evolution of the Fourier mode $ { |B_y(k_x, k_y)|}$; (c)  Time evolution of the Fourier mode $ { |B_z(k_x, k_y)|}$;  (d) Time evolution of the Fourier mode of density; $k_x = \frac{2\pi}{20}$, $k_y = \frac{2\pi}{20}$.}
%\label{fig:fire2D2fB}
%\end{figure}
%
%
%\begin{figure}[htbp]
%\center{
%\subfigure[]{\includegraphics[scale=0.49]{fnB_Bx_2D.png}}
%\subfigure[]{\includegraphics[scale=0.49]{fnB_By_2D.png}}
%\subfigure[]{\includegraphics[scale=0.49]{fnB_Bz_2D.png}}
%\subfigure[]{\includegraphics[scale=0.49]{fnB_density_2D.png}}
%}
%\caption{{\bf Mirror instability 2D, $\delta f$, $fn{\mathbf B}$ formulation, $1.5\times 10^5$ particles, $\alpha = 10^{-3}$ }: (a) Time evolution of the Fourier mode $ { |B_x(k_x, k_y)|}$; (b)  Time evolution of the Fourier mode $ { |B_y(k_x, k_y)|}$; (c)  Time evolution of the Fourier mode $ { |B_z(k_x, k_y)|}$;  (d) Time evolution of the Fourier mode of density; $k_x = \frac{2\pi}{20}$, $k_y = \frac{2\pi}{20}$.}
%\label{fig:fire2D2fnB}
%\end{figure}

\subsection{Parametric instability}
\begin{figure}[htbp]
\center{
{\includegraphics[scale=0.4]{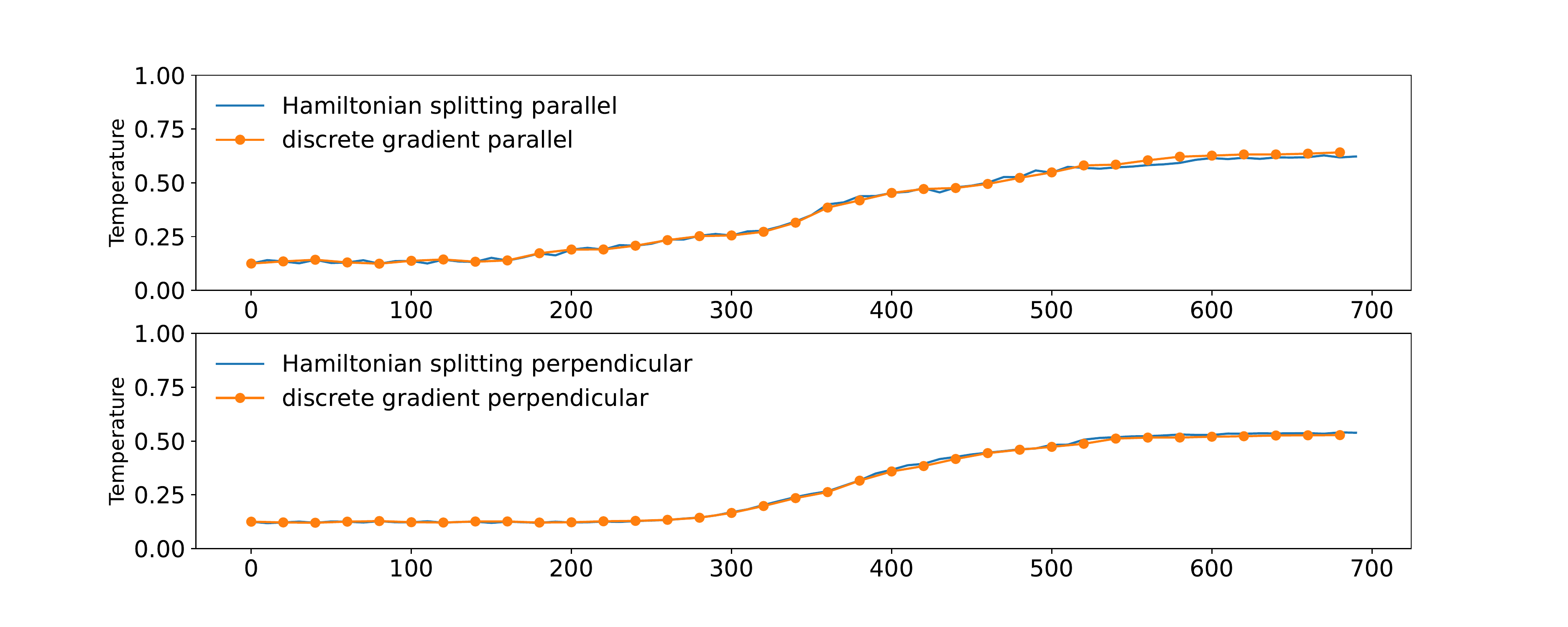}}
}
\caption{Time evolutions of parallel and perpendicular temperatures of ions given by schemes~\eqref{eq:LieI}-\eqref{eq:LieII}. }
\label{fig:para1D2}
\end{figure}

\begin{figure}[htbp]
\center{
{\includegraphics[scale=0.4]{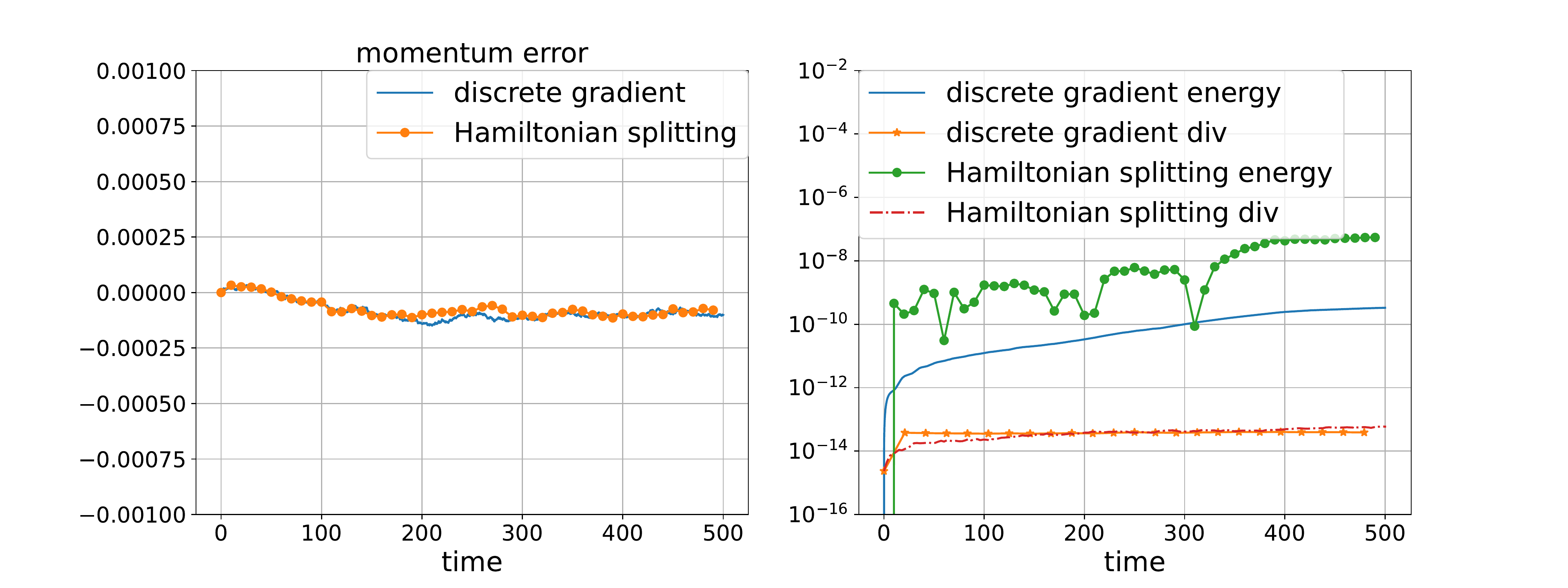}}
}
\caption{Time evolutions of momentum error (the first component) and relative energy error given by schemes~\eqref{eq:LieI}-\eqref{eq:LieII}. }
\label{fig:para1Derror}
\end{figure}

Finally we run the simulations of parametric instability~\cite{para} to compare and validate the schemes~\eqref{eq:LieI}-\eqref{eq:LieII}.
This instability produces a forward propagating acoustic wave and a backward Alfv\'en wave with a wave number smaller than that of the pump. 
We use the following initial conditions and numerical parameters to run simulations:
\begin{equation*}
\begin{aligned}
 &{\mathbf B} = (1, 0, 0) + \delta {\mathbf b}, \  \delta{\mathbf b} = (0, - \sin(0.196x), \cos(0.196x)), \ f  =  \frac{1}{\pi^{\frac{3}{2}} v_T^{\frac{3}{2}}} e^{-\frac{|{\mathbf v} - \delta{\mathbf u}|^2}{v_{T}^2} }, \, \kappa = 1,
\end{aligned}
\end{equation*}
where $v_T = 0.5$. 
Domain parameters are $[L_x,  L_y,  L_z] = [128, 1, 1]$, time step size is $\Delta t = 0.01$, total particle number is $10^6$, quadrature points in each cell are $[4, 2, 2]$, the degrees of basis polynomials are $[3, 1, 1]$, grids and degrees of shape functions are $[500, 3, 3]$,  $[2,2,2]$, respectively, cartesian mapping case is considered, and full $f$ method is used. 
As in~\cite{para}, $\delta {\mathbf u}$ is determined in the following way:
$$
\delta{\mathbf u} = -(\omega_0/k_0)\delta {\mathbf b}, \quad k_0^2 = \omega_0^2/(1-\omega_0), \quad k_0 = 0.196.
$$
By solving $\omega_0$ from $k_0^2 = \omega_0^2/(1-\omega_0)$, we get $\omega_0 = 0.17773094298487538$ (left-handed circularly polarized waves).  
Time evolutions of (parallel, perpendicular) temperatures of ions are plotted in Fig.~\ref{fig:para1D2}. We can see that instabilities happen during the simulations, the temperature behaviors of both schemes are close to the results in~\cite{para}. 
In Fig.~\ref{fig:para1Derror}, we present the momentum and relative energy error. 
The relative energy error of scheme~\eqref{eq:LieII} is at the level of $10^{-10}$, and relative energy error of scheme~\eqref{eq:LieI}   is at the level of $10^{-8}$. Momentum errors of both schemes~\eqref{eq:LieI}-\eqref{eq:LieII} are at the level of $10^{-4}$. Note that for scheme~\eqref{eq:LieII}, the degrees of shape functions are at least $[2,2,2]$ to guarantee the convergence of the discrete gradient methods of sub-step $xv$. As in the previous numerical tests, the weak errors of $\nabla \cdot {\mathbf B} = 0$ are at the level of iteration tolerance. 

\section{Conclusion}\label{sec:conclusion}
In this work, we propose two geometric methods for a hybrid model with kinetic ions and mass-less electrons based on a bracket structure. The conservation properties of energy, divergence free condition of magnetic field, positivity of the density, and quasi-neutrality relation are  investigated and compared for this two methods~\eqref{eq:LieI}-\eqref{eq:LieII}. Extensive numerical tests are done to validate and compare the schemes. Considering more realistic boundary conditions as~\cite{pml}, constructing momentum conserving schemes, and applying our numerical methods to more complicated physical simulations are future works.

\section*{Acknowledgements}
Simulations in this work were performed on Max Planck Computing \& Data Facility (MPCDF). Discussions with and inputs from Carlos Gonzalez, Dimitrios Kaltsas, Omar Maj,  Philip Morrison, and Anna Tenerani  are appreciated.

\end{document}